\documentclass[12pt,reqno]{amsart}
\usepackage[letterpaper,margin=1in]{geometry}
\usepackage[T1]{fontenc}
\usepackage{lmodern}
\usepackage{amsmath,amssymb,mathtools,mathrsfs,microtype,graphicx}
\usepackage{xcolor}
\usepackage[colorlinks=true,linkcolor=blue,citecolor=green!40!black,urlcolor=blue]{hyperref}
\hypersetup{pdftitle={Angular lenses and the failure of Simon's L2 conjecture},
 pdfauthor={Simon Becker and Mengxuan Yang}}
\newtheorem{theorem}{Theorem}[section]
\newtheorem{proposition}[theorem]{Proposition}
\newtheorem{lemma}[theorem]{Lemma}

\theoremstyle{remark}
\newtheorem{remark}[theorem]{Remark}
\numberwithin{equation}{section}
\newcommand{\R}{\mathbb R}
\newcommand{\C}{\mathbb C}
\newcommand{\Z}{\mathbb Z}
\newcommand{\dd}{\,d}
\newcommand{\1}{\mathbf1}
\newcommand{\supp}{\operatorname{supp}}
\newcommand{\dist}{\operatorname{dist}}
\newcommand{\rank}{\operatorname{rank}}

\author{Simon Becker}
\address{Bocconi University, Via Sarfatti 25, 20136 Milan, Italy}
\email{simon.becker@unibocconi.it}
\author{Mengxuan Yang}
\address{Texas A\&M University, College Station, TX 77843, USA}
\email{yangmx@tamu.edu}
\title[Angular lenses and pure point spectrum]{Lenses make many points:\\ Pure point spectrum disproves Simon's $L^2$-conjecture}
\begin{document}
      \begin{abstract}
    For every $d\ge2$ and $0\le\sigma<1/3$, we construct with LLM-assistance a smooth bounded nonnegative lens-type potential $V$ with arbitrarily small $\|V\|_\infty+\int_{\R^d}|x|^{1-d+\sigma}|V(x)|^2\dd x$, such that
$-\Delta+V$ has only pure point spectrum in $(0,\infty)$. Thus, we show that Barry Simon's $L^2$ conjecture, Problem~8 in \emph{Schr\"odinger operators in the twenty-first century}, is false.
	\end{abstract}
\maketitle
 \section{Introduction and main result}

For an integer $d\ge2$, let $\mathbb S^{d-1}$ be the unit sphere with surface measure $\dd\omega$. Given a real measurable function $V$ on $\R^d$ and $\sigma\ge0$, we
   define
\begin{equation*}
     \mathcal I_{d,\sigma}(V)
    \coloneqq\int_{\R^d}|x|^{1-d+\sigma}|V(x)|^2\dd x
  =\int_0^\infty r^\sigma
\int_{\mathbb S^{d-1}}|V(r\omega)|^2\dd\omega\dd r .
  \end{equation*}
   Simon's multidimensional $L^2$ conjecture claims that $\mathcal
  I_{d,0}(V)<\infty$ implies the presence of an absolutely continuous
 spectrum of infinite multiplicity at Lebesgue-almost every
positive energy for a self-adjoint realization of $-\Delta+V$. See \cite[Problem~8]{Simon00} and \cite[Conjecture~20.2]{Simon}. 

The conjecture is motivated by the hope of establishing general conditions under which a decaying potential preserves the absolutely continuous spectrum of the free Laplacian, see \cite[Section~20]{Simon}. Classical results of Kato, Kuroda, and Agmon~\cite{Kato59,Kuroda59,Agmon75} establish such preservation under the pointwise bound $|V(x)|\le C(1+|x|)^{-\beta}$ with $\beta>1$. In one dimension, Kiselev~\cite{Kiselev96} reached $\beta>3/4$, and Christ--Kiselev~\cite{CK98} and Remling~\cite{Remling98} subsequently reached $\beta>1/2$, whereas below $1/2$ absolutely continuous spectrum can be absent~\cite{Remling97}. The theorem of Deift--Killip~\cite{DK} then established preservation for every real $V\in L^2(\R)$, without a pointwise decay assumption. Together with the corresponding radial picture in higher dimensions, these results motivated the hope that the condition $\mathcal I_{d,0}(V)<\infty$ would provide a general theorem that implies the persistence of absolutely continuous spectrum. We refer to Safronov~\cite{Safronov08} and Denisov~\cite{Den} for a discussion of progress and partial results towards this result. 

 For a self-adjoint operator $H$, write $\sigma(H)$ and $\sigma_{\mathrm{ess}}(H)$ for its spectrum and essential spectrum. Let $P_{\mathrm{ac}}(H)$ and $P_{\mathrm{sc}}(H)$ be the projections onto its absolutely continuous and singular continuous
   subspaces, and define $P_c(H)\coloneqq P_{\mathrm{ac}}(H)+P_{\mathrm{sc}}(H)$.

   The following theorem disproves the conjecture and shows that even a stronger weighted $L^2$ condition does
not ensure any continuous, let alone absolutely continuous, spectral component on
      $[0,\infty)$.

  \begin{theorem}\label{thm:main}
      For $d\ge2$, $0\le\sigma<1/3$ and $\delta>0$,
  there is a smooth bounded nonnegative potential $V$ on $\R^d$ such that
\begin{equation}\label{eq:main-weight}
  \|V\|_\infty+\mathcal I_{d,0}(V)+\mathcal I_{d,\sigma}(V)<\delta
  \end{equation}
    and the operator $H_V\coloneqq-\Delta+V$ on $H^2(\R^d)$ satisfies
\begin{equation*}
     \sigma(H_V)=\sigma_{\mathrm{ess}}(H_V)=[0,\infty),\qquad
	 P_c(H_V)=0.
    \end{equation*}
	\end{theorem}

  \begin{remark}
 The continuous spectral component thus vanishes on $[0,\infty)$ even when the integrability condition includes the additional weight $|x|^\sigma$, $0<\sigma<1/3$. Pure point spectrum here means that the eigenfunctions form a complete orthonormal basis; the positive eigenvalues are dense in $[0,\infty)$.
  \end{remark}

\begin{figure}[tbp]
    \centering
\includegraphics[width=.7\textwidth]{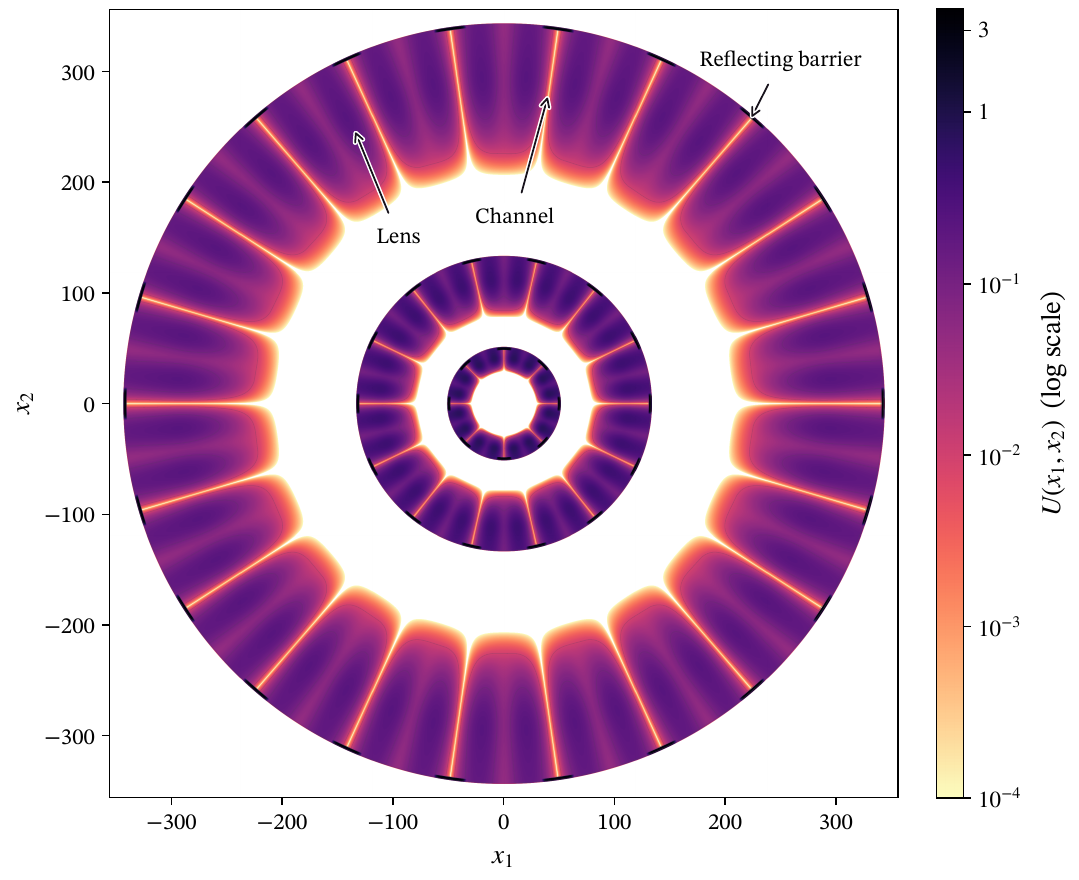}
\caption{Three illustrative annular blocks of our
      potential.} \label{fig:potential-blocks}
	\end{figure}

	The planar construction of our counterexample is a sum of potentials supported on
      widely separated annuli. See Figure~\ref{fig:potential-blocks}. Each annular
      potential consists of a weak angular lens and stronger reflecting barriers. The lens concentrates
    waves in narrow angular channels that meet the barriers. These names are illustrated in the figure.
The barriers consist of periodic pieces whose spectral gaps cover any prescribed bounded energy interval,
while their heights remain uniformly bounded. The energy intervals
 increase sufficiently slowly along the annular scales. In three dimensions, we extend the potential cylindrically and add disk-shaped barriers in the remaining direction. Tensor sums of the two- and three-dimensional operators give the construction in every dimension.

  For the two- and three-dimensional operators, we obtain pure point spectrum with finite eigenvalue multiplicities by applying the spectral averaging criterion of Combes--Hislop \cite[Theorem~3.2]{CH}, in the form stated
in \cite[Proposition~4.2]{FLM}. We derive the required
weighted resolvent bound from compactly supported approximate solutions. Proposition~\ref{prop:assembly} gives this step of the
    argument.

\medskip \noindent\textbf{Outline of the paper.}
Section~\ref{sec:spectral-averaging} derives a weighted resolvent bound from compactly supported approximate solutions and uses spectral averaging to obtain pure point spectrum.
Section~\ref{sec:planar} constructs the planar lenses and reflecting barriers, proves the required uniform estimates, and assembles the two-dimensional potential.
Section~\ref{sec:dimensions} establishes cylindrical confinement in three dimensions and completes the proof of Theorem~\ref{thm:main} in all dimensions by tensor sums and dilation.

\medskip \noindent\textbf{AI disclosure:} Our interest in this conjecture was renewed when it remained among the few of Simon's problems left untouched by OpenAI. This, together with the scarcity of positive results even under stronger restrictions on the potential, led us to suspect that the conjecture might
	be false. We developed the first construction of the potential in discussions with ChatGPT, obtaining a proof of absence of absolutely continuous spectrum in a bounded interval. The argument was rather different from the final manuscript and showed through some iteration procedure that the mass of spectral
   measures was concentrated on sets of increasingly small sets.

	We then conjectured that $[0,\infty)$ might nevertheless always be the union of the absolutely
	continuous and singular continuous spectra. When prompted with this question, ChatGPT showed that a
	rather intricate adjustment of the parameters of the potential, though not a change of potential, could
yield point spectrum in a bounded interval through a
   mechanism very similar to the one used in the final
  proof and already including the spectral-averaging argument. Very much to our
surprise, ChatGPT was then able to extend this result to the full spectrum by pushing the parameters of the construction further. We thoroughly rewrote the resulting argument and manuscript, and used ChatGPT again to help edit the final version. We take full responsibility for the mathematical content of this article.

\medskip \noindent\textbf{Notation.} Our inner products are
antilinear in the first entry. We use unnormalized Lebesgue measure
    on $\mathbb S^1\coloneqq\R/(2\pi\Z)$ and write $\langle
    x\rangle\coloneqq(1+|x|^2)^{1/2}$. For an open set
$\Omega\subset\R^d$, write $\iota_\Omega:L^2(\Omega)\to L^2(\R^d)$
  for the map $f\mapsto\iota_\Omega f$, where $\iota_\Omega f$ equals $f$ on $\Omega$ and zero outside of it. We identify $f$ with $\iota_\Omega f$ when working in $L^2(\R^d)$. Subscripts $j$ and $b$, or no subscript, indicate the domain specified in each statement. We use $B_r\coloneqq\{x:|x|<r\}$ for the
      centered Euclidean ball and $B_r(y)\coloneqq y+B_r$ for its translate. For planar coordinates, define $\mathbf e_1\coloneqq(1,0)$. We denote by $\1_J(H)$ the spectral projection of a self-adjoint $H$ on a Borel set $J$.

\section{Spectral averaging}\label{sec:spectral-averaging}

  The goal of the next proposition is to outline the strategy of the argument and prove pure point
spectrum by spectral averaging. We obtain the required weighted resolvent bound from compactly supported
      approximate solutions and then apply the criterion of Combes
   and Hislop \cite[Theorem~3.2]{CH}, in the form stated in
   \cite[Proposition~4.2]{FLM}.

\begin{proposition}\label{prop:assembly}
Let $d\ge2$, let $B\ge0$ be a bounded real potential on $\R^d$, and let $J$ be a bounded energy interval. Fix positive constants
$\alpha,\beta,Q$. Suppose $R_j\ge1$, $b_{j+1}=3R_j$, and, for fixed $p>1$ and $K$,
\[
10R_j\le R_{j+1}\le K R_j^p,\qquad b_j<R_j/4.
\]
  Let $V_j\ge0$ be uniformly bounded potentials supported in $\{R_j<|x|<9R_j/4\}$. Suppose $Z_j\subset J$ are measurable and
$\sum_j|Z_j|<\infty$. For every $j$ and every energy $\lambda\in J\setminus Z_j$, assume there is a bounded linear map
      $Q_j(\lambda):L^2(B_{b_j})\to H^2(\R^d)$ such that every $f\in L^2(B_{b_j})$ gives an approximate solution $u\coloneqq Q_j(\lambda)f$ satisfying $ \supp u\subset
      B_{9R_j/4}$
\begin{equation*}
    \begin{gathered}
\|(-\Delta+B+V_j-\lambda)u-\iota_jf\|_2
  \le C R_j^{-\beta}\|f\|_2.
     \end{gathered}
  \end{equation*}
    The same excluded energies $Z_j$ work for all $f$.
The constant $C$ is independent of $j,\lambda,f$. If
	\begin{equation*}
 \alpha>d,\qquad \beta+1>(p+2)\alpha,\qquad Q>p^2\alpha,
     \end{equation*}
	then, for almost every real $t$,
    \[
  -\Delta+B+\sum_jV_j+t\langle x\rangle^{-2Q}
      \]
has pure point spectral type on $J$, with finite eigenvalue multiplicities.
 \end{proposition}

      \begin{proof}
      We first compare two finite-volume resolvents and then pass to the whole space. The annuli supporting the $V_j$ are disjoint,
 so $\|\sum_jV_j\|_\infty\le\sup_j\|V_j\|_\infty$. The boundedness of $B$ and of this sum shows that
$H\coloneqq-\Delta+B+\sum_jV_j$ is self-adjoint on $H^2(\R^d)$. For each scale, let $\Omega_j\coloneqq
   B_{3R_j}$. The Dirichlet operators $H_j^0$ and $H_j$ on $\Omega_j$ have potentials $B+V_j$ and
$B+\sum_{i\le j}V_i$, respectively. Their common domain is $H^2(\Omega_j)\cap H_0^1(\Omega_j)$. The
superscript zero therefore indicates that all potentials $V_i$ with $i<j$ have been omitted.

  The finite-volume resolvents can be large near their eigenvalues. We remove small neighborhoods of these eigenvalues before estimating them. By nonnegativity and the min--max principle, each counting function is bounded by
that of the Dirichlet Laplacian on $B_{3R_j}$. Scaling this ball to $B_1$ and applying Weyl's law \cite[Theorem~3.3.1]{LMP} gives, for every fixed $E$,
    \[
\rank\1_{(-\infty,E]}(H_j^0)
  +\rank\1_{(-\infty,E]}(H_j)\le C_E R_j^d.
\]
      In particular, the sets
\[
	 \mathcal B_j\coloneqq\{\lambda\in J:
	 \dist(\lambda,\sigma(H_j^0)\cup\sigma(H_j))<R_j^{-\alpha}\}
     \]
 satisfy
\[
|\mathcal B_j|\le CR_j^{d-\alpha},\qquad
   \sum_j\bigl(|\mathcal B_j|+|Z_j|\bigr)<\infty.
	\]
Here $\alpha>d$ and $R_{j+1}\ge10R_j$ ensure summability of the first term. Borel--Cantelli now gives the full-measure set
 \[
      J_*\coloneqq J\setminus\limsup_{j\to\infty}(\mathcal B_j\cup Z_j),
\qquad |J\setminus J_*|=0.
  \]
For the rest of the proof, fix $\lambda\in J_*$. There is an index $j_0(\lambda)$ such that, for every $j\ge j_0(\lambda)$,
   \[
  \mathcal R_j^0\coloneqq(H_j^0-\lambda)^{-1},\qquad
\mathcal R_j\coloneqq(H_j-\lambda)^{-1},\qquad
    \|\mathcal R_j^0\|+\|\mathcal R_j\|\le2R_j^\alpha.
   \]
For every $j\ge j_0(\lambda)$, we also have $\lambda\notin Z_j$. The assumption therefore provides a map $Q_j(\lambda)$ satisfying the stated support and error bounds for every $f\in L^2(B_{b_j})$. The same $J_*$ and
$j_0(\lambda)$ work for every $f\in L^2(B_{b_j})$.

We next show that a solution with right-hand side
supported in $B_{b_j}$ is small near the boundary of $\Omega_j$. The relevant outer annulus and the error of the local inverse are
	\[
\mathcal O_j\coloneqq\{9R_j/4<|x|<3R_j\},\qquad
	 E_j\coloneqq(H_j^0-\lambda)Q_j-\iota_j.
  \]
	The functions $Q_jf$ belong to the Dirichlet domain because they vanish outside $B_{9R_j/4}$. Applying $\mathcal R_j^0$ to the definition of $E_j$ gives
	\[
    \mathcal R_j^0\iota_j=Q_j-\mathcal R_j^0E_j,
\qquad \|E_j\|\le CR_j^{-\beta}.
	\]
The first term on the right vanishes on $\mathcal O_j$. Consequently,
\[
      \|\1_{\mathcal O_j}\mathcal R_j^0\iota_j\|
    \le\|\mathcal R_j^0\|\|E_j\|
  \le CR_j^{\alpha-\beta}.
\]
    To restore the earlier annular potentials, we use
$U_{<j}\coloneqq\sum_{i<j}V_i$. Their supports and norms satisfy
\[
	 \supp U_{<j}\subset B_{b_j},\qquad
	 \|U_{<j}\|_\infty\le\sup_i\|V_i\|_\infty.
\]
  Thus, for every $f\in L^2(B_{b_j})$, the new right-hand side $g_{j,f}\coloneqq\iota_jf-U_{<j}\mathcal
 R_j\iota_jf$ is still supported in $B_{b_j}$. The resolvent identity gives
\[
     \mathcal R_j\iota_jf=\mathcal R_j^0g_{j,f},\qquad
 \|g_{j,f}\|\le(1+C R_j^\alpha)\|f\|.
  \]
	The preceding estimate for $\mathcal R_j^0$ therefore yields
 \[
\|\1_{\mathcal O_j}\mathcal R_j\iota_jf\|
	 \le CR_j^{\alpha-\beta}(1+C R_j^\alpha)\|f\|
      \le CR_j^{2\alpha-\beta}\|f\|.
    \]
 This estimate holds for every $f\in L^2(B_{b_j})$ with the same constant and the same excluded energies.

    We can now cut off the finite-volume solution inside $\Omega_j$. Let $\Theta_j$ be a smooth radial function with values in $[0,1]$, equal to one on $B_{5R_j/2}$ and zero outside $B_{11R_j/4}$, whose derivatives of order
$k$ are bounded by $C_kR_j^{-k}$. For $u\coloneqq\mathcal R_j\iota_jf$, the equation on $\mathcal O_j$ is $(-\Delta+B-\lambda)u=0$. The localized energy estimate \cite[Lemmas~3.5--3.6]{FLM} gives
	\[
    \|\nabla u\|_{\supp\nabla\Theta_j}
   \le C\bigl(R_j^{-1}+\sqrt{\max\{\lambda,0\}}\bigr)
  \|u\|_{\mathcal O_j}
      \le C_J\|u\|_{\mathcal O_j}.
\]
This follows by testing the equation with a cutoff times $\overline u$ and
      absorbing the term containing one derivative of that cutoff. Nonnegativity of
 $B$ is sufficient, and no derivative of $B$ occurs. The compactly supported
 approximate solution and its residual are
	\[
     v_jf\coloneqq\Theta_j\mathcal R_j\iota_jf,\qquad
    r_jf\coloneqq\iota_jf-(H-\lambda)v_jf.
     \]
 Since $\Theta_j=1$ on the ball and $H=H_j$ on its support,
      \[
\begin{aligned}
r_jf&=-[-\Delta,\Theta_j]u
   =2\nabla\Theta_j\cdot\nabla u+(\Delta\Theta_j)u,\\
    \|r_jf\|&\le CR_j^{-1}\|u\|_{\mathcal O_j}
	 \le CR_j^{2\alpha-\beta-1}\|f\|.
\end{aligned}
    \]

We next compare the solutions for the same $f$ at successive radii. Both $v_jf$ and its residual are supported in $B_{b_{j+1}}$,
	where $b_{j+1}=3R_j$. Moreover, $v_jf$ belongs to the next Dirichlet domain, $H_{j+1}v_jf=Hv_jf$, and $\Theta_{j+1}v_jf=v_jf$.
    Applying $\mathcal R_{j+1}$ to $(H_{j+1}-\lambda)v_jf=f-r_jf$, we obtain
\[
  \begin{aligned}
  v_jf&=\mathcal R_{j+1}(f-r_jf) \text{ and }
     v_{j+1}f-v_jf=\Theta_{j+1}\mathcal R_{j+1}r_jf.
   \end{aligned}
\]
   Here and below the same $f$ is extended by zero when viewed in a larger ball.
    The last identity bounds the change from one scale
      to the next:
      \[
    \|v_{j+1}f-v_jf\|
\le CR_{j+1}^\alpha R_j^{2\alpha-\beta-1}\|f\|
  \le CR_j^{(p+2)\alpha-\beta-1}\|f\|.
   \]
      By our assumption, the exponent is negative. Hence, this bound is summable in $j$. On each fixed space $L^2(B_{b_m})$, the
	operators $v_j$ therefore converge in norm as $j\to\infty$, with $j\ge\max(m,j_0(\lambda))$. Their limit is denoted by
      $G_{\lambda,m}$. For every $f$ in that space,
      \[
 \begin{gathered}
  v_jf\longrightarrow G_{\lambda,m}f,\qquad
 (H-\lambda)v_jf=\iota_mf-r_jf\longrightarrow\iota_mf,\\
	(H-\lambda)G_{\lambda,m}f=\iota_mf,\qquad
 \|G_{\lambda,m}\|\le R_{\max(m,j_0(\lambda))}^{\alpha}+C.
  \end{gathered}
     \]
The second line follows from closedness of $H$. Restricting the same
      sequence $v_j$ to a smaller ball also shows that these limits
 agree whenever both are applied to the same $f$.

	It remains to combine the operators $G_{\lambda,m}$ into a single bounded operator with a decaying weight on the right-hand side. We use
 \[
	 \mathcal C(x)\coloneqq\langle x\rangle^{-Q},\qquad
      b_0\coloneqq0,\qquad
  \mathcal A_m\coloneqq B_{b_m}\setminus B_{b_{m-1}}.
    \]
      The bound for $G_{\lambda,m}$ involves $R_m$, whereas the smallest radius on $\mathcal A_m$ is $b_{m-1}=3R_{m-2}$. The upper scale bound must therefore be used twice:
 \[
    R_m\le KR_{m-1}^p\le K^{p+1}R_{m-2}^{p^2}
\le Cb_{m-1}^{p^2}.
    \]
   Consequently, for all sufficiently large $m$, the coefficients
 $d_m\coloneqq\|G_{\lambda,m}\| \|\mathcal
    C\|_{L^\infty(\mathcal A_m)}$ satisfy
  \[
   d_m\le C_\lambda b_{m-1}^{p^2\alpha-Q},\qquad
 \sum_m d_m^2<\infty.
\]
The last conclusion uses $Q>p^2\alpha$ and geometric growth of the radii. The
finitely many earlier terms are finite as well. For any $f\in L^2(\R^d)$, the
series
\[
  F_\lambda f\coloneqq
 \sum_mG_{\lambda,m}(\mathcal C\1_{\mathcal A_m}f)
  \]
   converges absolutely in $L^2$, since
\[
\sum_m\|G_{\lambda,m}(\mathcal C\1_{\mathcal A_m}f)\|
	 \le\sum_m d_m\|\1_{\mathcal A_m}f\|
\le\left(\sum_m d_m^2\right)^{1/2}\|f\|.
  \]
    If $F_{\lambda,M}f$ denotes the sum through $m=M$, then
\[
   F_{\lambda,M}f\longrightarrow F_\lambda f,\qquad
(H-\lambda)F_{\lambda,M}f
	=\mathcal C\1_{B_{b_M}}f\longrightarrow\mathcal C f.
      \]
Closedness of $H$ now gives
    \[
	 F_\lambda f\in H^2(\R^d),\qquad
	 (H-\lambda)F_\lambda f=\mathcal C f.
	\]
 Thus $F_\lambda$ is bounded and satisfies
 $(H-\lambda)F_\lambda=\mathcal C$.

 This right inverse gives the resolvent bound needed for spectral averaging. For
   $\eta>0$, the function $m_\eta(s)\coloneqq s/(s-i\eta)$ has absolute value at
most one on $\R$. The spectral theorem and the preceding identity imply
\[
 \begin{gathered}
     (H-\lambda-i\eta)^{-1}\mathcal C
    =(H-\lambda-i\eta)^{-1}(H-\lambda)F_\lambda
  =m_\eta(H-\lambda)F_\lambda,\\
\sup_{\eta>0}\|(H-\lambda-i\eta)^{-1}\mathcal C\|
\le\|F_\lambda\|<\infty.
  \end{gathered}
	\]
  This holds for every $\lambda$ in the full-measure set $J_*$. No uniform bound as
      $\lambda$ varies is required.

  The remaining hypotheses of the spectral averaging criterion are compactness and dense range. Because
    $\mathcal C$ is strictly positive,
   \[
  \bigl(\operatorname{Ran}\mathcal C\bigr)^\perp
=\ker\mathcal C=\{0\}.
  \]
  Its range is therefore dense. For $z\notin\R$, boundedness of the potential gives $(H-z)^{-1}:L^2\to H^2$. If $\chi_L\coloneqq\1_{B_L}$, Rellich compactness shows that $\chi_L\mathcal C(H-z)^{-1}\mathcal C$ is compact. The
full operator is its norm limit, because
	\[
  \|\mathcal C(H-z)^{-1}\mathcal C
-\chi_L\mathcal C(H-z)^{-1}\mathcal C\|
  \le\frac{\langle L\rangle^{-Q}}{|\Im z|}\longrightarrow0.
   \]
All three assumptions of \cite[Proposition~4.2]{FLM} now hold with $h_0\coloneqq H$ and $h_1\coloneqq\mathcal C^2$. That criterion
gives pure point spectral type on $J$, with finite eigenvalue multiplicities, for almost every coupling in
  \[
   h_0+th_1=H+t\langle x\rangle^{-2Q}.
 \]
   This is the claimed operator.
  \end{proof}

   \section{Uniform planar estimates}
  \label{sec:planar}

   We fix $0\le\sigma<1/3$. For functions $f$ supported in a disk
      $B_b$, with $b>1$ large, we will place the next annular
	potential at radius $R$. We choose the parameters
    \begin{equation}\label{eq:planar-parameters}
	 \begin{gathered}
\ell\coloneqq\frac{1-3\sigma}{16},\qquad r_*\coloneqq\frac32,
\qquad\eta\coloneqq2\ell,\\
	A\coloneqq\lceil b^2\rceil,\quad R\coloneqq A^{r_*},\quad
M\coloneqq\lceil A^{(r_*+\ell)/2}\rceil,\quad d_A\coloneqq A^\ell.
   \end{gathered}
 \end{equation}
Thus $b^3\le R\le2b^3$ for large $b$. The parameters balance the costs of the lens and barrier, as verified in \eqref{eq:planar-cost}. Fix an upper energy bound $E>0$. All $o(1)$ terms in this section tend to zero as
$b\to\infty$, uniformly over the energy and angular ranges stated in each estimate. Constants may depend on $E$ unless stated otherwise.

   Choose a fixed even $\chi\in C_c^\infty(( -\pi/2,\pi/2);[0,1])$, equal to one on $[-\pi/3,\pi/3]$, and extend it periodically. Set
\[
	 W_a(\phi)\coloneqq
2a\bigl(1-\chi(\phi)\operatorname{sech}^2(\sqrt a\,\phi)\bigr),
\qquad a\ge0,
      \]
      where the formula is read on $[-\pi,\pi]$ and extended periodically. The smooth parameter $a(t)$, defined below with the
  bounds \eqref{eq:planar-schedule}, is zero near $t\le1$, equals $A$ on a terminal interval of length comparable to $A^{-1}$ before $t=2$, and is extended as $A$ for $t\ge2$. It satisfies
   \[
    0\le a(t)\le A,\qquad
 \int_1^2a(t)^2\,dt\le CA.
\]
We choose an even $F_A\in C^\infty(\mathbb S^1;[0,1])$, equal to one on $|\phi|\le
      A^{-1/2+\eta}$ and zero on $|\phi|\ge2A^{-1/2+\eta}$.
Lemma~\ref{lem:radial-reflector} supplies a smooth reflecting potential
	\[
	q_{\rm bar}=q_{A,E}\in C_c^\infty((0,d_A);[0,3]).
 \]
It consists of finitely many periodic pieces, chosen
for the energy interval $[0,E]$, with smooth transitions. We also fix $e_{d_A}\in C^\infty(\R;[0,1])$ satisfying
      \[
e_{d_A}(x)=1\ \text{for }x\le d_A+1,\qquad
  e_{d_A}(x)=0\ \text{for }x\ge d_A+2.
     \]
 We construct $e_{d_A}$ by shifting a fixed smooth cutoff by $d_A$. Hence, for
every integer $j\ge0$, $\|\partial_x^j q_{\rm bar}\|_\infty+ \|\partial_x^j
	e_{d_A}\|_\infty\le C_{E,j}$, where $C_{E,j}$ is independent of $d_A$. In polar
   coordinates $(r,\vartheta)$, the potential is the sum of a lens and a reflecting
barrier. With $x\coloneqq r-2R$, these are
\[
  V_{\rm lens}\coloneqq e_{d_A}(x)\frac{M^2}{r^2}W_{a(r/R)}(M\vartheta),\qquad
     V_{\rm bar}\coloneqq e_{d_A}(x)q_{\rm bar}(x)F_A(M\vartheta),\qquad
V_{b,E}\coloneqq V_{\rm lens}+V_{\rm bar}.
\]
    We abbreviate $V_{b,E}$ to $V_b$ while $E$ is fixed, and extend $V_b$ by zero near the origin. Periodicity and the integer value of $M$ give
\[
      \int_0^{2\pi}|F_A(M\vartheta)|^2\,d\vartheta
	 =\frac1M\int_0^{2\pi M}|F_A(\phi)|^2\,d\phi
   =\int_0^{2\pi}|F_A(\phi)|^2\,d\phi\le CA^{-1/2+\eta}.
  \]
  Since $\|W_a\|_\infty\le2a$, the two contributions satisfy
    \[
\begin{gathered}
 \mathcal I_{2,\sigma}(V_{\rm lens})
    \le C_\sigma M^4R^{\sigma-3}\int_1^2a(t)^2\,dt
	+C_\sigma M^4A^2(d_A+2)R^{\sigma-4}
\le C_\sigma M^4A R^{\sigma-3},\\
  \mathcal I_{2,\sigma}(V_{\rm bar})
  \le C_\sigma d_AR^\sigma A^{-1/2+\eta}.
\end{gathered}
\]
   The first bound uses the integral estimate for $a(t)$ and $A(d_A+2)/R=O(A^{\ell-1/2})$. The second uses radial thickness $d_A$ and the preceding angular integral. Thus, for sufficiently large $b$, $V_b$ takes values in $[0,4]$ and $\operatorname{supp}V_b\subset\{R<r<9R/4\}$, and
\begin{equation}\label{eq:planar-cost}
	\begin{aligned}
\mathcal I_{2,\sigma}(V_b)
 &\le C_\sigma\left(M^4A R^{\sigma-3}
	+d_AR^\sigma A^{-1/2+\eta}\right)\\
	&\le C_\sigma\left(A^{2\ell+1-r_*(1-\sigma)}
+A^{r_*\sigma+\ell+\eta-1/2}\right)
	 \le C_\sigma A^{-5\ell}.
  \end{aligned}
\end{equation}
The constants in this cost bound are independent of $E$. The exponents satisfy
 \[
\begin{gathered}
    2\ell+1-r_*(1-\sigma)=-6\ell<0,\qquad
r_*\sigma+\ell+\eta-\tfrac12=-5\ell<0,\\
   \frac{M^2A}{R^2}\le CA^{\ell+1-r_*}
     =CA^{\ell-1/2}\longrightarrow0.
\end{gathered}
  \]

	The next proposition supplies the approximate inverses needed in Proposition~\ref{prop:assembly}. Given $f\in L^2(B_b)$, it constructs $u\coloneqq Q_b(\mu)f$ that is zero outside $B_{9R/4}$ and solves $(-\Delta+V_b-\mu)u=f$ up to an arbitrarily small error.

	\begin{proposition}\label{prop:planar-block}
  Fix $E>0$ and $B,T>0$. For all sufficiently large $b$ there is a measurable $Z_b\subset[R^{-1/4},E]$ with $|Z_b|\le C_E R^{-T}$
 such that, for every $\mu\in[R^{-1/4},E]\setminus Z_b$, a bounded linear operator
    \[
	 Q_b(\mu):L^2(B_b)\longrightarrow H^2(\mathbb R^2)
\]
	has, for every $f\in L^2(B_b)$, the properties
	\begin{equation*}
	 \supp Q_b(\mu)f\subset B_{9R/4},\qquad
  \|(-\Delta+V_{b,E}-\mu)Q_b(\mu)f-\iota_bf\|_2\le R^{-B}\|f\|_2.
  \end{equation*}
For each fixed $E$, the same potential $V_{b,E}$ works for all $B,T$. The lower bound on $b$, the map $Q_b$, and the set $Z_b$ may depend on these parameters. The same excluded energies $Z_b$ work for every $f$. The height, support, and cost bounds for $V_{b,E}$ are independent of $E$, once $b$ is
sufficiently large.
      \end{proposition}

The lens concentrates waves into narrow angular regions where the barrier reflects them. We first show that the angular ground state is concentrated in these regions.
   Lemma~\ref{lem:lens-profiles} constructs waves whose angular profiles follow that state through the lens. At its outer edge, Lemma~\ref{lem:reflected-waves} joins them to solutions decaying in a suitable periodic piece of the barrier. A cutoff then makes them vanish beyond the barrier with a small residual.
The reflected waves are expressed in terms of free solutions inside the annulus to match them to the source.

    For every $f\in L^2(B_b)$, \S\ref{sec:planar-matching} combines these waves to define $Q_b(\mu)f$, supported in $B_{9R/4}$, and proves the support and error bounds of Proposition~\ref{prop:planar-block}.
The excluded set $Z_b$ is common to all sources, and the estimates for each periodic angular profile include every Fourier mode.

\subsection{Angular confinement}\label{sec:angular-confinement}
	We need a gap above the angular ground state that grows with $a$. We also need
concentration of that state when $a=A$. For $a\ge0$ and $\zeta\in\R$, we set
$D\coloneqq-i\partial_\phi$ and
    \[
    K_{a,\zeta}\coloneqq(D+\zeta)^2+W_a\quad\hbox{on }H^2(\mathbb S^1).
     \]
We write $\lambda_j(K_{a,\zeta})$, $j\ge1$, for its eigenvalues in nondecreasing order, counted with multiplicity, and set
	$\Lambda(a,\zeta)\coloneqq\lambda_1(K_{a,\zeta})$. The real
      parameter $\zeta$ shifts each angular Fourier frequency
 $k\in\Z$ to $k+\zeta$. The potential
   $W_a=2a(1-\chi\operatorname{sech}^2(\sqrt a\,\phi))$ is smooth at
	$a=0$, since $\operatorname{sech}^2$ is even and analytic near zero.
    Its derivatives satisfy
	\[
	(1+a)^j\|\partial_a^jW_a\|_\infty\le C_j(1+a)
	\quad\text{for every integer }j\ge0.
    \]
  For $a\ge1$, this follows by differentiating in $y\coloneqq\sqrt a\,\phi$: every polynomial times a derivative of $\operatorname{sech}^2y$ is bounded. For $0\le a\le1$, it follows from smoothness on a compact parameter
interval.

We first estimate the gap for $\zeta=0$. On the line, set
    \[
     \begin{gathered}
H_a^{\rm line}\coloneqq-\partial_\phi^2+2a\tanh^2(\sqrt a\,\phi),
    \qquad \mathsf A_a\coloneqq\partial_\phi+\sqrt a\tanh(\sqrt a\,\phi),\\
\mathsf A_a^*\mathsf A_a=H_a^{\rm line}-a,
\qquad\mathsf A_a\mathsf A_a^*=-\partial_\phi^2+a.
 \end{gathered}
\]
For $a>0$, the kernel of $\mathsf A_a$ is spanned by the normalized function $\psi_a(\phi)\coloneqq
   c_a\operatorname{sech}(\sqrt a\,\phi)$. The polar decomposition of $\mathsf A_a$ and the lower bound
      $\mathsf A_a\mathsf A_a^*\ge a$ therefore give
\[
     \langle u,H_a^{\rm line}u\rangle
   \ge a\|u\|^2+a\|u-\langle\psi_a,u\rangle\psi_a\|^2.
   \]
Choose a smooth partition $\xi_0^2+\xi_1^2=1$ on the circle, with $\xi_0=1$ on $[-\pi/2,\pi/2]$ and
      $\xi_0=0$ for $|\phi|\ge3\pi/4$. The function $\xi_0u$ extends by zero to the line, where $W_a\ge2a\tanh^2(\sqrt a\,\phi)$. On the support of $\xi_1$, we have $W_a=2a$. If $u\perp\xi_0\psi_a$, then $\xi_0u\perp\psi_a$
 on the line. The partition identity gives
  \[
   \begin{aligned}
\langle u,K_{a,0}u\rangle
=\sum_{j=0}^1\langle\xi_ju,K_{a,0}\xi_ju\rangle
  -\sum_{j=0}^1\|\xi_j'u\|^2 \ge2a\bigl(\|\xi_0u\|^2+\|\xi_1u\|^2\bigr)-C\|u\|^2
     =(2a-C)\|u\|^2.
\end{aligned}
	\]
    The min--max principle thus gives $\lambda_2(K_{a,0})\ge2a-C$.

   Choose a fixed cutoff $\chi_0\in C_c^\infty((-\pi/3,\pi/3))$, equal to one near zero. On its support the line and circle potentials agree, so $\chi_0\psi_a$ is a periodic quasimode at energy $a$:
    \[
\|(K_{a,0}-a)\chi_0\psi_a\|\le Ce^{-c\sqrt a},\qquad
\|\chi_0\psi_a\|=1+O(e^{-c\sqrt a}).
     \]
The residual is the cutoff commutator, supported where $\psi_a$ and its derivative are exponentially small. The second-eigenvalue bound identifies the nearby eigenvalue
   as the ground eigenvalue for large $a$. Let $g_a$ be the positive normalized ground state at $\zeta=0$. We have proved
\[
\Lambda(a,0)=a+O(e^{-c\sqrt a}),\qquad
	 \lambda_2(K_{a,0})-\Lambda(a,0)\ge ca
      \quad\text{for large }a.
\]
      On a compact interval of $a$, the positive periodic ground
  state is simple, so its gap has a positive minimum. The gap is therefore at least $c(1+a)$ for every $a\ge0$.

The diamagnetic inequality \cite[\S9, (9.24)--(9.27)]{SimonDia} and the trial state $g_a$ give
	\[
   |\partial_\phi|u||\le|(D+\zeta)u|,\qquad
  \Lambda(a,0)\le\Lambda(a,\zeta)
\le\langle g_a,K_{a,\zeta}g_a\rangle=\Lambda(a,0)+\zeta^2.
      \]
      The form inequality
      \[
K_{a,\zeta}\ge(1-\tau)K_{a,0}+\zeta^2-\zeta^2/\tau
	 \qquad(0<\tau<1)
      \]
    with a fixed small $\tau$ preserves a gap of size $ca$ for large
 $a$ and $|\zeta|\le1/4$. Continuity of the eigenvalues on the
remaining compact interval gives constants $g_*>0$ and $\zeta_0\in(0,1/4)$ such that
    \begin{equation}\label{eq:planar-cell-gap}
\lambda_2(K_{a,\zeta})-\Lambda(a,\zeta)\ge g_*(1+a),
   \qquad |\zeta|\le\zeta_0,\quad a\ge0.
\end{equation}
  We also have $0\le\Lambda(a,\zeta)\le C(1+a)$ in this range. For later use, the diamagnetic inequality gives
\begin{equation}\label{eq:planar-angular-lower}
  K_{a,\zeta}\ge ca
	\quad\text{for every real }\zeta\text{ and all sufficiently large }a.
   \end{equation}

  We choose the radial schedule as follows. Let $s_0\in
C^\infty([1,2];[0,1])$ equal one near $1$, equal $2-t$ near $2$, and
   satisfy $c(2-t)\le s_0(t)\le C(2-t)$. Choose a smooth
 $F:[0,\infty)\to[0,\infty)$ such that $F(u)=0$ for $u\le1$, $F(u)=u$
  for $u\ge2$, and $0\le F(u)\le u$. Set
      \[
  h_A\coloneqq(1+A)^{-1},\quad
  s_A(t)\coloneqq h_A+(1-h_A)h_AF(s_0(t)/h_A),\quad
v(t)\coloneqq s_A(t)^{-1},\quad a(t)\coloneqq v(t)-1.
   \]
 Then $a=0$ near $t=1$ and $a=A$ on a terminal interval of length comparable to $A^{-1}$. We extend $a$ by these constant values
   outside $[1,2]$. Since $s_A\asymp h_A+2-t$, direct differentiation and integration give, for $j\ge1$ and $q>1$,
    \begin{equation}\label{eq:planar-schedule}
\begin{gathered}
      |\partial_t^jv|\le C_jv^{j+1},\qquad
 \int_1^2v(t)\,dt\le C\log(2+A),\qquad
\int_1^2v(t)^q\,dt\le C_q(1+A)^{q-1}.
\end{gathered}
     \end{equation}
	For the derivative bound, on the moving transition interval $s_A\asymp h_A$ and $|s_A^{(j)}|\le C_jh_A^{1-j}$. Away from it, all derivatives are bounded. Thus
  $|s_A^{(j)}|\le C_js_A^{1-j}$ everywhere, which gives the bound for $v=1/s_A$. In particular, $\int_1^2a(t)^2\,dt\le CA$.

We fix $\zeta$ satisfying \eqref{eq:planar-cell-gap}
      and follow the ground state as the lens parameter $a(t)$ varies. For $1\le t\le2$, we write $\Lambda(t)\coloneqq\Lambda(a(t),\zeta)$ and $K(t)\coloneqq K_{a(t),\zeta}$. We write $\Pi(t)$ for the orthogonal projection
onto the ground state and define its reduced inverse by
   \[
T(t)\coloneqq(K(t)-\Lambda(t)+\Pi(t))^{-1}-\Pi(t).
\]
   It inverts $K-\Lambda$ on the orthogonal complement of the ground state:
    \[
 T=(1-\Pi)T(1-\Pi),\qquad (K-\Lambda)T=1-\Pi.
	\]
Primes in the following estimates denote derivatives with respect to $t$. Since
      $a(1)=0$ and $|\zeta|\le\zeta_0<1/2$, the ground state at $t=1$ is the
  normalized constant function. We fix its phase for all $t$ by solving
\[
	 g'=[\Pi',\Pi]g,\qquad g(1)=(2\pi)^{-1/2}.
\]
	This choice gives $\Pi g=g$, $\|g\|=1$, and $\langle g,g'\rangle=0$. For the
      symmetry $\mathsf Jz(\phi)\coloneqq\overline{z(-\phi)}$, we have $\mathsf JK=K\mathsf J$ and hence $\mathsf J\Pi=\Pi\mathsf
J$. Thus $\mathsf Jg$ satisfies the same initial-value problem:
    \[
(\mathsf Jg)'=[\Pi',\Pi]\mathsf Jg,\qquad \mathsf Jg(1)=g(1).
	 \]
Uniqueness gives $\mathsf Jg=g$. The isolated-eigenvalue perturbation theory and parallel transport of
    \cite[Chapters~II, \S4, and VII, \S1]{Kato} apply on the fixed domain $H^2(\mathbb S^1)$. Here
  \[
	 \begin{gathered}
K'=a'\partial_aW_a,\qquad
\|\partial_t^jK\|\le C_jv^{j+1}\quad\text{for }j\ge1,\qquad
   \|T\|\le(g_*v)^{-1},\\
  \Lambda'=\langle g,K'g\rangle,\qquad g'=-TK'g.
      \end{gathered}
	 \]
    To estimate derivatives, we fix $t_0$ and use the contour
   $\Gamma\coloneqq\{z:|z-\Lambda(t_0)|=g_*v(t_0)/8\}$. For
      $|t-t_0|\le c/v(t_0)$, the scales $v(t)$ and $v(t_0)$ are
    comparable. With $\mathcal R(z,t)\coloneqq(K(t)-z)^{-1}$, we have
   \[
      \|K(t)-K(t_0)\|\le g_*v(t_0)/16,\quad
\|\mathcal R\|_{L^2\to L^2}\le C/v(t_0),\quad
 \partial_t\mathcal R=-\mathcal R K'\mathcal R.
   \]
For fixed $A$, inversion on the fixed domain gives
 \[
   \begin{gathered}
	 \Pi=-\frac1{2\pi i}\int_\Gamma\mathcal R\,dz,\qquad
  \mathcal R,\Pi,T\in C^\infty\bigl(\mathcal B(L^2,H^2)\bigr).
\end{gathered}
 \]
      Repeated differentiation, and $g'=[\Pi',\Pi]g$, give
\begin{equation}\label{eq:planar-cell-derivatives}
	\begin{gathered}
\|\partial_t^j\Pi\|+\|\partial_t^jg\|\le C_jv^j,
	 \qquad\|\partial_t^jT\|\le C_jv^{j-1},\qquad
	 |\partial_t^j\Lambda|\le C_jv^{j+1}
 \quad\text{for }j\ge0.
      \end{gathered}
\end{equation}
 For the reduced inverse, these bounds follow by differentiating
$T=(K-\Lambda+v\Pi)^{-1}-v^{-1}\Pi$: the inverse has norm $O(v^{-1})$, and the
  derivative of order $j\ge1$ of the inverted operator has norm $O(v^{j+1})$. The
      constants are uniform in $t_0,A,\zeta$ under \eqref{eq:planar-cell-gap}.
   Elliptic regularity also makes $g$ smooth in $\phi,t$.

	To prove concentration at $a=A$, write $\widetilde
g\coloneqq\chi_0e^{-i\zeta\phi}\psi_A$, with the cutoff
used above. It is periodic and
\[
  \|(K_{A,\zeta}-A)\widetilde g\|\le Ce^{-c\sqrt A},\qquad
\|\widetilde g\|=1+O(e^{-c\sqrt A}).
\]
	Indeed, the residual is $[(D+\zeta)^2,\chi_0]e^{-i\zeta\phi}\psi_A$, supported where $\psi_A$ is
 exponentially small. The gap identifies the nearby eigenvalue as the ground eigenvalue, so $\Lambda(A,\zeta)=A+O(e^{-c\sqrt A})$. With $c_A\coloneqq\langle
  g(2),\widetilde g\rangle$, it also gives
\[
\|\widetilde g-c_Ag(2)\|\le Ce^{-c\sqrt A},\quad |c_A|\ge\tfrac12,
   \qquad
\int_{|\phi|\ge A^{-1/2+\eta}}|\psi_A(\phi)|^2\,d\phi
  \le Ce^{-cA^\eta}.
\]
   The triangle inequality now yields $\|(1-F_A)g(2)\|\le
Ce^{-cA^\eta}$. Here the norm is taken in $L^2(\mathbb S^1)$, with
    all Fourier modes included.

\subsection{Approximate solutions through the lenses}\label{sec:lens-solutions}
We now construct oscillatory solutions whose angular part follows the varying ground state $g(t)$ from \S\ref{sec:angular-confinement}. The construction is often called adiabatic. We will control its error by a small
	parameter and keep every constant uniform in the angular cutoff.

      For the error exponents $B,T$ in Proposition~\ref{prop:planar-block}, we set $L\coloneqq T+3$ and choose an integer $P\ge1$ such that
	\begin{equation}\label{eq:planar-expansion-order}
\frac{\ell P}{2}>r_*(B+L+3),\qquad r_*=\frac32.
	\end{equation}
The order $P$ is independent of $b$, $\mu$, and the angular channel. We fix
  $\mu\in[R^{-1/4},E]$. With positive constants $C_0,C_1$ to be chosen in
    \S\ref{sec:planar-matching}, we define
\[
	 N\coloneqq\lceil C_0b+C_1\log R\rceil,\qquad n\coloneqq2N+1.
\]
We use the Fourier basis $e_m(\vartheta)\coloneqq(2\pi)^{-1/2}e^{im\vartheta}$, indexed by $m\in\Z$. For $|m|\le N$, we define $\zeta\coloneqq m/M$ and $\mathcal E_mz(\vartheta)\coloneqq e^{im\vartheta}z(M\vartheta)$. Since
   $M>2N$ for large $b$, the sets of Fourier modes $m+M\mathbb Z$ are disjoint and
 \begin{equation}\label{eq:planar-ladder-orthogonality}
      \left\|\sum_{|m|\le N}\mathcal E_mz_m\right\|^2
     =\sum_{|m|\le N}\|z_m\|^2.
\end{equation}
      With $D_\vartheta\coloneqq-i\partial_\vartheta$, we have
$D_\vartheta^2\mathcal E_m=M^2\mathcal E_m(D+\zeta)^2$. Each function $z_m$ may contain every Fourier
   mode on $\mathbb S^1$. The choice of scales in \eqref{eq:planar-parameters} gives
\[
\begin{gathered}
N\le Cb\le CA^{1/2},\qquad M^2\ge A^{r_*+\ell},\\
     N^2/M^2\le CA^{-1/2-\ell},\qquad N^2/R\le CA^{-1/2}.
	\end{gathered}
 \]
 These bounds ensure \eqref{eq:planar-cell-gap} and the uniform free-wave estimates
      below.

      The polar unitary $u\mapsto r^{1/2}u$ changes the Laplacian to
      $-\partial_r^2+r^{-2}(D_\vartheta^2-1/4)$. On the interval $t\coloneqq
	r/R\in[1,2]$, we set $\kappa\coloneqq\sqrt\mu$ and define the radial momentum and phase by
	\[
k_m(t)\coloneqq\left(\mu-\frac{M^2\Lambda(t)}{R^2t^2}
    +\frac1{4R^2t^2}\right)^{1/2},\qquad
      p_m\coloneqq k_m/\kappa,\qquad S_m(t)\coloneqq\int_1^t k_m(u)\,du.
	\]
To show that the radial momentum $k_m$ stays close to $\sqrt\mu$, we set $\delta\coloneqq M^2/(R^2\mu)$
and $\nu_m\coloneqq p_m^2-1$. For the entire range $R^{-1/4}\le\mu\le E$, we have
      \[
     \begin{gathered}
\nu_m=-\delta t^{-2}\Lambda+(4R^2\mu t^2)^{-1},\qquad
   \delta(1+A)\le CA^{\ell-1/8}=o(1),\qquad |\nu_m|=o(1),\\
|\nu_m^{(j)}|\le C_j\delta v^{j+1}+C_jR^{-7/4}
  \le C_jv^j\quad\text{for }j\ge1.
      \end{gathered}
	 \]
	Differentiating $p_m^{\pm1}=(1+\nu_m)^{\pm1/2}$ gives $|\partial_t^kp_m^{\pm1}|+|\partial_t^kp_m^{\pm1/2}|\le C_kv^k$. The estimate for $\nu_m$ keeps $p_m$ bounded away from zero. We also define $\epsilon$, which is
independent of $t$, by
   \[
\begin{gathered}
      \epsilon\coloneqq R\kappa/M^2,\qquad \epsilon^2\delta=M^{-2},\\
   \epsilon\log(2+A)\le C_EA^{-\ell}\log(2+A)=o(1),\qquad
	\epsilon\delta v=\frac{v}{R\kappa}\le CA^{-5/16}=o(1).
\end{gathered}
\]
	For $s\in\{+1,-1\}$,
  \[
 e^{-siRS_m}\partial_t^2(e^{siRS_m}q)
=q''+2siRk_mq'+siRk_m'q-R^2k_m^2q.
   \]
      Conjugation by $e^{siRS_m}$ and multiplication by $R^2/M^2$ give
   \begin{equation*}
 \mathcal F_s(q)\coloneqq t^{-2}(K-\Lambda)q
      -si\epsilon(2p_mq'+p_m'q)-\epsilon^2\delta q''.
     \end{equation*}
   In particular, we retain the term $-\epsilon^2\delta q''(t)$. It
      comes from the radial derivative
   $\partial_r^2q(r/R)=R^{-2}q''(t)$.

 We follow the isolated-band method of \cite[Theorem~2.1 and \S4.5]{PST} and \cite[Theorems~1.1 and 1.3]{MS}. The weighted derivative estimates needed for this family are
	proved below. For a fixed angular mode $m$, we abbreviate $p\coloneqq p_m$ and define
	\[
\begin{gathered}
\Pi^\perp\coloneqq1-\Pi,\qquad \mathscr L\coloneqq t^{-2}(K-\Lambda),\qquad
\mathscr G\coloneqq1+v^{-1}\mathscr L,\\
	\mathscr D_p\coloneqq p^{1/2}\partial_t^2(p^{-1/2}\,\cdot)
      =\partial_t^2-\frac{p'}p\partial_t
 +\frac{3(p')^2}{4p^2}-\frac{p''}{2p}.
 \end{gathered}
\]
  Since $\Lambda$ is the ground eigenvalue, $\mathscr G\ge1$. Since $2p(p^{-1/2}w)'+p'p^{-1/2}w=2p^{1/2}w'$, the change $q=p^{-1/2}w$ gives
 \[
     p^{1/2}\mathcal F_s(p^{-1/2}w)
  =\mathscr Lw-2si\epsilon p\,w'
 -\epsilon^2\delta\mathscr D_pw.
     \]
For $k\ge1$, the bounded remainder $\mathscr R_k$ in
   \[
\partial_t^k\mathscr L
   =(-1)^k(k+1)!t^{-k}\mathscr L+\mathscr R_k,
	\qquad \|\mathscr R_k\|\le C_kv^{k+1},
  \]
is controlled by \eqref{eq:planar-cell-derivatives}: the unbounded part comes only from differentiating $t^{-2}$, since $\partial_t^j(K-\Lambda)$ is bounded for $j\ge1$. For $k\ge1$, differentiating $\mathscr Lg=0$ and $\mathscr LT=t^{-2}\Pi^\perp$ gives
  \[
  \begin{aligned}
\mathscr L\partial_t^kg
=-\sum_{h=1}^k\binom kh
	(\partial_t^h\mathscr L)\partial_t^{k-h}g, \qquad 
   \mathscr L\partial_t^kT
=\partial_t^k(t^{-2}\Pi^\perp)
	-\sum_{h=1}^k\binom kh
(\partial_t^h\mathscr L)\partial_t^{k-h}T.
	\end{aligned}
\]
At $k=0$, $\mathscr Gg=g$ and $\mathscr
GT=T+v^{-1}t^{-2}\Pi^\perp$. The formula for $\partial_t^h\mathscr L$ reduces each sum to lower-order graph bounds and bounded
      remainders. Induction, using \eqref{eq:planar-cell-derivatives}, therefore gives
    \[
\|\mathscr G\partial_t^kg\|\le C_kv^k,\qquad
\|\mathscr G\partial_t^kT\|_{L^2\to L^2}\le C_kv^{k-1}
 \quad\text{for }k\ge0.
     \]
      Every constant is independent of $m$ and $N$. We need bounds that also control the angular derivatives of the error. For a profile $u(t)\in H^2(\mathbb S^1)$,
	smooth in $t$, we therefore use the weighted norm
   \[
 \|u\|_{k,A;\mathscr G}
\coloneqq\max_{0\le h\le k}\sup_{1\le t\le2}
 v(t)^{-h}\|\mathscr G(t)\partial_t^hu(t)\|.
    \]
The next lemma estimates the angular profiles' distance from the ground state and their residual.
   The scalar form near the endpoints permits matching of values and radial derivatives to the solutions on both sides of the lens.

   \begin{lemma}\label{lem:lens-profiles}
 With the parameters and operators defined above, for
	sufficiently large $b$ there are smooth profiles
$q_s(t,\phi)$, $s\in\{+1,-1\}$, such that, for every integer $k\ge0$,
	\[
 \begin{gathered}
\|q_s-p_m^{-1/2}g\|_{k,A;\mathscr G}
      \le C_{P,k}\epsilon\log(2+A),\\
      \|\mathscr G(t)\mathcal F_s(q_s)(t)\|
  \le C_P\epsilon v(t)[\epsilon\log(2+A)]^P
     \quad\text{for }1\le t\le2.
\end{gathered}
    \]
   Near $t=1$ and $t=2$, these profiles have the form $q_s=B_sg$,
  where $B_s(t)\coloneqq\langle g(t),q_s(t)\rangle$ and
  \[
	 B_s(1)=p_m(1)^{-1/2},\qquad
      B_s(2)=p_m(2)^{-1/2}+O_P(\epsilon\log(2+A)),\qquad
	 B_s'=O_P(A).
     \]
   Moreover, $q_-=\mathsf Jq_+$ and $B_-=\overline{B_+}$. All
constants are uniform in $A,\mu,m,N$ in the ranges fixed above.
	\end{lemma}

   \begin{proof}
   We start with $g$ and correct the equation separately in the directions orthogonal and parallel to $g$. The reduced inverse $T$ acts on the orthogonal component. For
      each sign $s$, we define
 \begin{gather}
 \label{eq:planar-scalar-recurrence}
 \mathcal A_su(t)\coloneqq(\mathcal V_su)(t)+g(t)\beta_u(t), \qquad 
 \mathcal V_su \coloneqq\epsilon t^2T
\left(2si\,p\,u'+\epsilon\delta\mathscr D_pu\right),\\
   \beta_u(t) \coloneqq\int_1^t
 \left(\langle g'(\tau),(\mathcal V_su)(\tau)\rangle
+\frac{si\epsilon\delta}{2p(\tau)}
	\langle g(\tau),(\mathscr D_pu)(\tau)\rangle\right)\,d\tau.  
    \end{gather}
    The integral determines the coefficient of $g(t)$ and vanishes at $t=1$. The weights $v(t)^{-h}$ in this norm account for the derivative growth proved in \eqref{eq:planar-cell-derivatives}. Since scalar multiplication
commutes with $\mathscr G$, the product rule gives, pointwise in $t$,
\[
  \|\mathscr G\partial_t^h u'\|
      \le v^{h+1}\|u\|_{h+1,A;\mathscr G},\qquad
 \|\mathscr G\partial_t^h(\mathscr D_pu)\|
      \le C_hv^{h+2}\|u\|_{h+2,A;\mathscr G}.
      \]
    We estimate the two terms in $\mathcal A_su$ using the bounds for $\mathscr G\partial_t^kg$ and $\mathscr G\partial_t^kT$ just proved. For $h\ge1$, differentiation gives
   \[
    \beta_u^{(h)}
  =\partial_t^{h-1}\left(\langle g',\mathcal V_su\rangle+
\frac{si\epsilon\delta}{2p}\langle g,\mathscr D_pu\rangle\right).
   \]
  We first use the bounds for $T$ to estimate the component orthogonal to $g$:
\[
\|\mathscr G\partial_t^h(\mathcal V_su)(t)\|
\le C_k\epsilon v(t)^h(1+\epsilon\delta v(t))
\|u\|_{k+2,A;\mathscr G}\quad(0\le h\le k).
	\]
A derivative of order $j$ on $T$ contributes $v^{j-1}$, while the remaining derivatives of $pu'$ contribute $v^{h-j+1}$. The powers cancel. The second derivative term leaves one extra factor of $v$, accompanied by $\epsilon\delta$. The formula for $\beta_u$, the bounds for $g$, and \eqref{eq:planar-schedule}
then give
  \[
	\begin{gathered}
\|\beta_u\|_\infty
\le C\epsilon\bigl(\log(2+A)+\epsilon\delta(1+A)
+\delta(1+A)\bigr)\|u\|_{2,A;\mathscr G},\\
    v^{-h}|\beta_u^{(h)}|
  \le C_k\epsilon(1+\epsilon\delta v+\delta v)
	 \|u\|_{k+2,A;\mathscr G} \qquad \text{for } 1\le h\le k,\\
  \|g\beta_u\|_{k,A;\mathscr G}
\le C_k\max_{h\le k}\sup_{1\le t\le2}v(t)^{-h}|\beta_u^{(h)}(t)|.
	 \end{gathered}
\]
   For the first bound, the two integrals that occur are
	$\int_1^2v\,dt=O(\log(2+A))$ and $\int_1^2v^2\,dt=O(A)$. For the
    second, differentiating the integrand $h-1$ times gives $C_k\epsilon v^h(1+\epsilon\delta v)+C_k\epsilon\delta v^{h+1}$. Since
$\delta(1+A)=o(1)$ and $\epsilon=o(1)$, combining these estimates for $\mathcal
      A_su=\mathcal V_su+g\beta_u$ yields
      \[
     \|\mathcal A_su\|_{k,A;\mathscr G}
  \le C_k\epsilon\log(2+A)\|u\|_{k+2,A;\mathscr G}.
	 \]
 We apply this correction $P$ times, where $P$ was fixed in \eqref{eq:planar-expansion-order}. We define
      \[
	\begin{gathered}
z_0\coloneqq g,\qquad z_j\coloneqq\mathcal A_sz_{j-1},\qquad
w_j\coloneqq\sum_{h=0}^jz_h\quad\text{for }0\le j\le P,\\
     q_s\coloneqq p^{-1/2}w_P,\qquad B_s\coloneqq\langle g,q_s\rangle.
  \end{gathered}
\]
For each $0\le j\le P$, induction gives
   \[
\|z_j\|_{k,A;\mathscr G}
  \le C_{j,k}[\epsilon\log(2+A)]^j\|g\|_{k+2j,A;\mathscr G}
	 \le C_{j,k}[\epsilon\log(2+A)]^j.
    \]
    Thus only derivatives through order $k+2P$ are needed, and
\[
      \|q_s-p^{-1/2}g\|_{k,A;\mathscr G}
    \le C_{P,k}\epsilon\log(2+A).
  \]
Since $\langle g,\mathcal V_su\rangle=0$ and $\langle g,g'\rangle=0$, differentiating
\eqref{eq:planar-scalar-recurrence} gives
    \[
\langle g,(\mathcal A_su)'\rangle
  =-\langle g',\mathcal V_su\rangle+\langle g',\mathcal V_su\rangle
   +\frac{si\epsilon\delta}{2p}\langle g,\mathscr D_pu\rangle
	=\frac{si\epsilon\delta}{2p}\langle g,\mathscr D_pu\rangle.
\]
      Summing the defining equations yields
 \[
    \Pi^\perp w_P=\mathcal V_sw_{P-1},\qquad
 \Pi w_P'=\frac{si\epsilon\delta}{2p}
   \Pi\mathscr D_pw_{P-1}.
\]
	In particular,
   \[
\mathscr Lw_P=2si\epsilon p\,\Pi^\perp w_{P-1}'
   +\epsilon^2\delta \Pi^\perp\mathscr D_pw_{P-1}.
   \]
Substitution in the normalized equation, with $w_P-w_{P-1}=z_P$,
   gives the exact remainder
  \[
\begin{aligned}
  p^{1/2}\mathcal F_s(q_s)
    &=-2si\epsilon p\,\Pi^\perp z_P'-\epsilon^2\delta\mathscr D_pz_P,\\
\|\mathscr G\mathcal F_s(q_s)(t)\|
   &\le C_P\epsilon v(t)[\epsilon\log(2+A)]^P
  (1+\epsilon\delta v(t))\\
   &\le C_P\epsilon v(t)[\epsilon\log(2+A)]^P.
\end{aligned}
  \]
 Here $\mathscr G$ commutes with $\Pi^\perp$ and with scalar functions of $t$. These estimates include every Fourier mode on
$\mathbb S^1$. Their constants are uniform in $A,\mu,m,N$. Since $\|W_a\|_\infty\le CA$, elliptic regularity also gives
 $\|\partial_t^kq_s\|_{H^2}\le C_{P,k}A^{k+1}$. Undoing the operator scaling and $dr=R\,dt$ multiplies the residual norm by
$(M^2/R^2)\sqrt R=M^2/R^{3/2}$. Using $\int_1^2v^2\,dt\le C(1+A)$ and $\epsilon\log(2+A)\le C_EA^{-\ell/2}$ for large $A$, we obtain
\begin{equation*}
\begin{gathered}
      U_{m,s}(r,\vartheta)\coloneqq\mathcal E_m
     \bigl(e^{siRS_m(r/R)}q_s(r/R)\bigr),\\
\frac{M^2}{R^{3/2}}\epsilon(1+A)^{1/2}
	[\epsilon\log(2+A)]^P
=\kappa R^{-1/2}(1+A)^{1/2}[\epsilon\log(2+A)]^P
    \le C_{P,E}\kappa A^{-1/4-\ell P/2},\\
	\|(-\partial_r^2+r^{-2}(D_\vartheta^2-1/4)+V_b-\mu)U_{m,s}\|
 _{L^2((R,2R)\times\mathbb S^1,\,dr\,d\vartheta)}
\le C_{P,E}\kappa A^{-1/4-\ell P/2}.
	\end{gathered}
      \end{equation*}

Near $t=1$ and $t=2$, the function $a(t)$ is constant by construction. Thus $K$, $\Pi$, and $g$ are constant there, and $g'=0$. For
      a scalar function $a_0(t)$,
	\[
     \partial_t(a_0g)=a_0'g,\quad \mathscr D_p(a_0g)=(\mathscr D_pa_0)g,
	 \quad Tg=0\quad\Longrightarrow\quad\mathcal V_s(a_0g)=0.
     \]
    It follows inductively that $z_j$ is a scalar multiple of $g$ on
these intervals. At $t=1$, the integral defining $\beta_u$ is zero,
so $z_j(1)=0$ for $j\ge1$. Near $t=2$, the integral may be nonzero.
	For the amplitudes $B_s=\langle g,q_s\rangle$ defined above, we
	obtain
\[
\begin{gathered}
   q_s=B_sg,\qquad B_s(1)=p_m(1)^{-1/2},\qquad
 B_s(2)=p_m(2)^{-1/2}+O_P(\epsilon\log(2+A)),\qquad B_s'=O_P(A).
      \end{gathered}
      \]
The scalar amplitudes can still vary, since $p_m$
varies. The operator $\mathsf J$ commutes with $K,T$
  and changes $i$ to $-i$. Applying it to the recurrence therefore gives $q_-=\mathsf Jq_+$ and $B_-=\overline{B_+}$.
	\end{proof}

\subsection{Reflection at the outer barrier}\label{sec:barrier-reflection}
   We now join the two lens solutions to a solution that decays inside the reflecting barrier. Recall the
 shifted radial coordinate $x=r-2R$. For $r\ge2R$, the angular parameter is $a=A$. We construct the
 barrier from a finite family of periodic potentials whose spectral gaps cover $[0,E]$. Their heights are
  bounded independently of $E$.

      \begin{lemma}[A finite family of periodic gaps]\label{lem:periodic-gaps}
      There are finitely many smooth real periodic functions $q_0,\ldots,q_{n_E}$, periods $T_i>0$, and
 compact intervals $J_i\subset[0,E]$ whose union is
$[0,E]$, with the following properties. We may take
   $q_0=3$ and $J_0=[0,\min\{2,E\}]$. For $i\ge1$,
$0\le q_i\le1$. For every $i$ and every $\nu\in J_i$,
  the system
      \[
     Y'=\begin{pmatrix}0&1\\q_i(x)-\nu&0\end{pmatrix}Y
	\]
has an exponential dichotomy. If $U_i(x,y;\nu)$ is its propagator and $P_{i,s}(x;\nu),P_{i,u}(x;\nu)$ are
    its stable and unstable projections, then, for constants $K_E,\gamma_E>0$,
    \[
     \begin{aligned}
    \|U_i(x,y;\nu)P_{i,s}(y;\nu)\|
     &\le K_E e^{-\gamma_E(x-y)},&&x\ge y,\\
  \|U_i(x,y;\nu)P_{i,u}(y;\nu)\|
&\le K_E e^{-\gamma_E(y-x)},&&x\le y.
    \end{aligned}
\]
 The phase of each $q_i$ can be chosen so that a unit unstable
 eigenvector $v_{i,u}(\nu)$ of $U_i(T_i,0;\nu)$ satisfies
  \[
   |(v_{i,u}(\nu))_1|\ge c_E>0\qquad(\nu\in J_i).
      \]
	For $q_0=3$, the same conclusions hold on every compact subset of $(-\infty,3)$; its use at arbitrarily negative energies will
   instead be controlled directly by an energy estimate.
	\end{lemma}

   \begin{proof}
   The constant potential $3$ has characteristic exponents $\pm\sqrt{3-\nu}$ for $\nu<3$, so it handles $J_0$. Fix an energy $\nu_0>1$ and set
\[
  k_0=\sqrt{\nu_0},\qquad k_1=\sqrt{\nu_0-1},\qquad
    a=\frac\pi{2k_0},\qquad b=\frac\pi{2k_1}.
\]
Consider first the $(a+b)$-periodic potential equal to
zero on the first interval of length $a$ and to one on
the following interval of length $b$. At energy
	$\nu_0$, its monodromy in the coordinates
      $Y=(u,u')^{\mathsf T}$ is
 \[
	\begin{pmatrix}0&1/k_1\\-k_1&0\end{pmatrix}
  \begin{pmatrix}0&1/k_0\\-k_0&0\end{pmatrix}
  =\begin{pmatrix}-k_0/k_1&0\\0&-k_1/k_0\end{pmatrix}.
 \]
      The two multipliers are real, reciprocal, and
have moduli strictly above and below one. The unstable eigenline is the first
 coordinate axis.

   Periodic smoothing by a nonnegative mollifier leaves
   the potential between zero and one and changes it by
      an arbitrarily small amount in $L^1$ over one
 period. The variation-of-constants formula and
 Gronwall's inequality show that its monodromy
    converges to the displayed matrix. Hence a
sufficiently small smoothing preserves strict hyperbolicity and the nonzero first
coordinate of the unstable line. The same properties hold on an open energy
	neighborhood of $\nu_0$, by continuous dependence of the monodromy on the energy.
Shrinking this neighborhood to a compact interval gives uniform bounds on its multipliers, eigenprojections, and propagators within one period. Iteration by whole periods gives the stated exponential estimates.

 Apply this construction at the points of the compact interval $[2,E]$ when $E>2$, and take a finite subcover by smaller compact intervals contained in the corresponding open neighborhoods. We may intersect the resulting intervals with $[0,E]$. Together with $J_0$ they cover $[0,E]$. Taking the minimum and
	maximum of the finitely many positive constants proves the lemma.
	\end{proof}

   We record explicitly the stability of the decaying solution under the small centrifugal term that occurs after separation of the angular variable. Fix a member of the preceding family and an energy in its compact interval.
  Abbreviate its propagator by $U$, its projections by $P_s,P_u$, and common dichotomy constants by $K,\gamma$. If $e:[0,\infty)\to\R$ is continuous and
    $\|e\|_\infty\le\delta$, write 
\[
B(x)=\begin{pmatrix}0&0\\e(x)&0\end{pmatrix}.
	 \]
  Choose a real unit vector $v_s\in\operatorname{ran}P_s(0)$. On the Banach space with norm
\[
  \|Y\|_\gamma=\sup_{x\ge0}e^{\gamma x/2}|Y(x)|,
\]
	consider the integral equation
   \begin{equation}\label{eq:reflector-fixed-point}
\begin{split}
  Y(x)={}&U(x,0)v_s
 +\int_0^x U(x,t)P_s(t)B(t)Y(t)\,\dd t -\int_x^\infty U(x,t)P_u(t)B(t)Y(t)\,\dd t.
      \end{split}
	\end{equation}
The norms of the two integral operators are at most $2K\delta/\gamma$ and $2K\delta/(3\gamma)$, respectively. Thus, for a sufficiently small $\delta_E>0$, their sum has norm at most $1/2$, uniformly over the finite family and
its energy intervals. The equation has a unique real solution, with
	\begin{equation}\label{eq:reflector-stable-solution}
    |Y(x)|\le C_Ee^{-\gamma_E x/2},\qquad
      c_E\le |Y(0)|\le C_E.
\end{equation}
 For the lower bound, apply $P_s(0)$ to \eqref{eq:reflector-fixed-point} at zero: the unstable integral disappears and
      $P_s(0)Y(0)=v_s$. Differentiating the convergent integrals proves that $Y$ solves
   \[
Y'=\begin{pmatrix}0&1\\q_i(x)+e(x)-\nu&0\end{pmatrix}Y.
\]
The estimates use only $\|e\|_\infty\le\delta_E$ and do not require a bound on its derivatives.

	\begin{lemma}\label{lem:radial-reflector}
Fix $E>0$ and $\ell>0$. Let $A\to\infty$, let $R=R(A)>0$, and let $c_A\ge0$ satisfy $c_A/R^2\to0$. Set $d_A=A^\ell$. For sufficiently large $A$ there is a function
    \[
     q_A\in C_c^\infty((0,d_A);[0,3])
  \]
with the following property. For every $\mu\in[0,E]$ and every real
      $c$ with $|c|\le c_A$, there is a real smooth function
	$w=w_{A,\mu,c}$ on $[0,\infty)$, vanishing near $d_A$ and beyond,
  such that
\begin{equation}\label{eq:radial-reflector-bounds}
\begin{gathered}
|w(0)|^2+|w'(0)|^2=1,\qquad
   \sup_{x\ge0}(|w(x)|+|w'(x)|)\le C_E e^{C_Ed_A},\\
	\left\|\left(-\partial_x^2+q_A(x)+\frac{c}{(2R+x)^2}-\mu\right)w
	\right\|_{L^2(0,\infty)}
 \le C_E e^{-c_Ed_A}.
 \end{gathered}
      \end{equation}
  The potential $q_A$ is independent of $\mu$ and $c$. For every fixed integer $j\ge0$, it may be chosen with $\|q_A^{(j)}\|_\infty\le C_{E,j}$.
   \end{lemma}

    \begin{proof}
  Use the finite family of Lemma~\ref{lem:periodic-gaps}. We concatenate one long piece of each $q_i$, in the order $i=0,\ldots,n_E$, with smooth transitions of fixed length two. The first piece starts at $x=2$, after a
smooth transition from zero, and the last is followed by a smooth transition to zero. Convex interpolation by fixed smooth cutoffs keeps the potential in $[0,3]$ and bounds every fixed number of its derivatives independently
      of the lengths.

We now choose these lengths. For the scalar first-order equation in this proof, the coefficient matrix has norm at most a fixed $C_0(E)$ for all sufficiently large $A$. Let $\gamma_E$ be as in \eqref{eq:reflector-stable-solution}. Choose positive numbers $a_0,\ldots,a_{n_E}$ successively so large that
 \[
\frac{\gamma_Ea_i}{4}
-C_0(E)\sum_{h<i}a_h\ge\frac{\gamma_Ea_i}{8}
\qquad(0\le i\le n_E).
\]
   For example, a geometric sequence with ratio larger than $1+8C_0(E)/\gamma_E$
suffices. Take 
    \[
  l_i=\frac{a_i}{2\sum_{h=0}^{n_E}a_h},\qquad L_i=l_i d_A,
   \qquad s_i=2+\sum_{h<i}(L_h+2).
  \]
The periodic part of the $i$th piece occupies $[s_i,s_i+L_i]$. Their
total length is $d_A/2$, and the finitely many transitions fit inside
 $(0,d_A)$ for large $A$. Moreover, for some $c_0(E)>0$ and all large
    $A$,
 \begin{equation}\label{eq:reflector-lengths}
\frac{\gamma_E L_i}{4}-C_0(E)s_i\ge c_0(E)d_A
	\qquad(0\le i\le n_E).
    \end{equation}
	We denote the resulting smooth concatenation by $q_A$.

   Given $\mu\in[0,E]$, choose $i$ with $\mu\in J_i$. Beginning at
 $x=s_i$, extend the $i$th periodic potential to the whole right
 half-line and perturb it there by $c/(2R+x)^2$. Its sup norm is at
      most $c_A/(4R^2)=o(1)$. The fixed-point construction
	\eqref{eq:reflector-fixed-point} therefore gives a real exact
solution whose Cauchy vector at $s_i$ has norm between two positive
constants and whose Cauchy vector for $x\ge s_i$ is bounded by $C_Ee^{-\gamma_E(x-s_i)/2}$. Continue this solution backwards from
$s_i$ to zero through the actual concatenated potential. Gronwall's inequality for the propagator and its inverse gives
      \[
|Y(0)|\ge c_Ee^{-C_0(E)s_i}.
      \]
Divide by $|Y(0)|$. The normalized solution, still
 denoted by $Y$, has $|Y(0)|=1$ and satisfies
    \[
   \begin{aligned}
|Y(x)|&\le e^{C_0(E)x},&&0\le x\le s_i,\\
|Y(x)|&\le C_Ee^{C_0(E)s_i}
	e^{-\gamma_E(x-s_i)/2},&&x\ge s_i
\end{aligned}
\]
	as long as the selected periodic equation is used on
   the right. In particular,
 \eqref{eq:reflector-lengths} makes its value and
derivative at $x\ge s_i+L_i/2$ at most
$C_Ee^{-c_Ed_A}$.

   Choose a fixed-shape smooth cutoff equal to one up
      to $s_i+L_i/2$ and zero after $s_i+L_i/2+1$.
These two points lie strictly inside the $i$th periodic
piece for all large $A$. Multiply the scalar component
  of $Y$ by this cutoff and call the result $w$. Where
      the cutoff is one, $w$ solves the exact equation
   with the actual $q_A$. Its residual is the
commutator $-\chi''Y_1-2\chi'Y_2$, supported on an interval of length one and bounded by $C_Ee^{-c_Ed_A}$. Its Cauchy data at zero remain normalized. The preceding estimates, also on the cutoff interval, prove
\eqref{eq:radial-reflector-bounds}.
 \end{proof}

 For the angular channel $m$, define
      \[
    c_m\coloneqq M^2\Lambda(A,m/M)-\tfrac14.
 \]
Since $\Lambda(A,m/M)\le CA$ and $M^2A/R^2=o(1)$, the centrifugal coefficient satisfies $|c_m|/R^2=o(1)$ uniformly in $|m|\le N$. Apply Lemma~\ref{lem:radial-reflector}
with $c_A=CM^2(A+1)$ and $q_{\rm bar}=q_A$, and denote the resulting scalar function by $w_m$. It has normalized Cauchy data at $x=0$, is zero before $x=d_A$, and
    satisfies the exact scalar equation except on its cutoff interval.

	We use the transformed radial operator $\mathcal L_b\coloneqq-\partial_r^2+r^{-2}(D_\vartheta^2-1/4)+V_b$. On the support of $w_m$, we have $a=A$ and $e_{d_A}=1$. The
actual barrier contains the angular factor $F_A$. Using the concentration estimate from \S\ref{sec:angular-confinement}, we compute
      \[
	\begin{aligned}
(\mathcal L_b-\mu)(w_m\mathcal E_mg(2))
&=\left(-w_m''+
    \left(q_{\rm bar}+\frac{c_m}{(2R+x)^2}-\mu\right)w_m\right)
    \mathcal E_mg(2)\\
  &\quad+q_{\rm bar}w_m\mathcal E_m((F_A-1)g(2)).
	\end{aligned}
	 \]
The first term is exponentially small by \eqref{eq:radial-reflector-bounds}. The concentration estimate gives $\|(F_A-1)g(2)\|\le Ce^{-cA^{2\ell}}$. Although $w_m$ may grow before reaching the appropriate periodic piece, its growth is bounded by $C_Ee^{C_Ed_A}$. Consequently the resulting error is bounded by
 \begin{equation}\label{eq:planar-barrier-error}
  \varepsilon_{{\rm bar},E}(A)\coloneqq
    C_Ee^{-c_EA^\ell}
      +C_EA^{\ell/2}\exp(C_EA^\ell-cA^{2\ell}).
	 \end{equation}
This is smaller than every negative power of $A$, and hence of $R$, for each fixed $E$. The angular term $M^2r^{-2}W_A$ remains in the scalar equation through its ground eigenvalue in $c_m$.

We use the amplitudes $B_s$ from \S\ref{sec:lens-solutions} to define two scalar waves near the endpoints:
\[
     \begin{gathered}
L_{m,s}(r)\coloneqq e^{siR(S_m(r/R)-S_m(2))}B_s(r/R),\qquad
\mathsf C_m(r)\coloneqq
  \begin{pmatrix}L_{m,+}(r)&L_{m,-}(r)\\L_{m,+}'(r)&L_{m,-}'(r)\end{pmatrix}.
\end{gathered}
	\]
The columns of $\mathsf C_m$ record the values and $r$-derivatives of these two waves. Thus $\mathsf C_m(r)(a,b)^{\mathsf T}$ gives the value and derivative of
    $aL_{m,+}+bL_{m,-}$ at $r$. Since $B_-=\overline{B_+}$, the estimates at $t=2$ give
\[
\begin{aligned}
\det\mathsf C_m(2R)
&=-2i\kappa p_m(2)|B_+(2)|^2
+\frac{B_+(2)\overline{B_+'(2)}-\overline{B_+(2)}B_+'(2)}R\\
  &=-2i\kappa\left(1+O_P\!\left(\epsilon\log(2+A)+
    \frac{A}{\kappa R}\right)\right).
	\end{aligned}
  \]
 The error tends to zero uniformly, since $A/(\kappa R)\le
    A^{-5/16}=o(1)$. Thus
\[
\begin{gathered}
\|\mathsf C_m(2R)\|\le C,\quad
\|\mathsf C_m(2R)^{-1}\|\le C/\kappa,\\
     \binom{\gamma_m}{\overline{\gamma_m}}
\coloneqq\mathsf C_m(2R)^{-1}\binom{w_m(0)}{w_m'(0)},\qquad
    c\le|\gamma_m|\le C/\kappa.
     \end{gathered}
      \]

    \subsection{Free radial waves and reflection}\label{sec:free-reflection}
We next match the lens solutions to free waves at the inner radius $r=R$. This requires estimates uniform over the growing range $|m|\le N$. We use $H_j^{(1)},H_j^{(2)}$ for the Hankel functions of the first and second kind, and $J_j$ the Bessel function of the first kind, all of order $j$. We normalize the
free waves as in \cite[\S10.17]{DLMF}:
   \[
     \alpha_m\coloneqq\pi|m|/2+\pi/4,\qquad
F_m^+(r)\coloneqq\sqrt{\pi\kappa r/2}\,e^{i\alpha_m}
	H_{|m|}^{(1)}(\kappa r),\qquad
	 F_m^-\coloneqq\overline{F_m^+}.
     \]
   We arrange their values and derivatives in the fundamental matrix
\[
\mathsf F_m(r)\coloneqq
\begin{pmatrix}F_m^+(r)&F_m^-(r)\\(F_m^+)'(r)&(F_m^-)'(r)\end{pmatrix},
      \qquad\det\mathsf F_m(r)=-2i\kappa.
  \]
The next lemma supplies the reflected waves used in the final proof. At $r=R$, each
wave is the sum of an outgoing free wave and an incoming free wave of the same
      amplitude. Beyond the barrier it vanishes. We need both a small error in the
equation and a lower bound for its mass to control the matching coefficients in
  \S\ref{sec:planar-matching}.

\begin{lemma}\label{lem:reflected-waves}
 For sufficiently large $b$, every $\mu\in[R^{-1/4},E]$ and $|m|\le N$ admit a complex number $\rho_m$ with $|\rho_m|=1$ and a function $\Psi_m$ on $[R,\infty)\times\mathbb S^1$ with the following properties. The function
	vanishes for $r\ge9R/4$, has angular Fourier support in $m+M\Z$, and satisfies
\[
 \binom{\Psi_m(R,\cdot)}{\partial_r\Psi_m(R,\cdot)}
	=\binom{F_m^+(R)+\rho_mF_m^-(R)}
  {(F_m^+)'(R)+\rho_m(F_m^-)'(R)}e_m.
\]
      For the transformed operator $\mathcal L_b$ defined in \S\ref{sec:barrier-reflection}, the residual
$r_m\coloneqq(\mathcal L_b-\mu)\Psi_m$ satisfies
 \[
\|r_m\|\le\delta_P\coloneqq
     C_{P,E}\kappa R^{-1/2}A^{1/2}[\epsilon\log(2+A)]^P
	+\varepsilon_{{\rm bar},E}(A).
\]
Furthermore,
	\[
	c_E R\le\int_R^{2R}\|\Psi_m(r,\cdot)\|_{L^2(\mathbb S^1)}^2\,dr
	\le C_E R.
   \]
    The last norm and the residual norm use the measure $dr\,d\vartheta$. The constants are uniform in $\mu,m,N$ in the stated ranges.
\end{lemma}

   \begin{proof}
We use the Hankel remainder estimate \cite[(10.17.13)--(10.17.15)]{DLMF} with truncation order one. For every integer $n\ge0$ and $x>0$, it gives
      \[
      \left|\sqrt{\pi x/2}\,e^{-i(x-\pi n/2-\pi/4)}H_n^{(1)}(x)-1\right|
    \le\frac{|n^2-1/4|}{x}
  \exp\!\left(\frac{|n^2-1/4|}{x}\right).
\]
	Apply this bound at $n=|m|$ and $n+1$, and use the derivative recurrence
      \cite[(10.6.2)]{DLMF},
      \[
(H_n^{(1)})'(x)=\frac{n}{x}H_n^{(1)}(x)-H_{n+1}^{(1)}(x).
\]
    If $x\ge1+|n^2-1/4|$, the exponential factors for both orders are bounded by absolute constants, since $|(n+1)^2-1/4|\le3(1+|n^2-1/4|)$. Consequently, for $\kappa r\ge1+|m^2-1/4|$,
  \[
F_m^+(r)=e^{i\kappa r}\left(1+O\!\left(\frac{1+m^2}{\kappa r}\right)\right),\qquad
 \kappa^{-1}(F_m^+)'(r)
  =e^{i\kappa r}\left(i+O\!\left(\frac{1+m^2}{\kappa r}\right)\right).
\]
The constants are independent of $m$. The condition on $\kappa r$ holds uniformly for $r\ge R/2$ and $|m|\le N$ when $b$ is large,
      because
\[
  \frac{N^2+1}{\kappa R}\le C_E R^{2/3-7/8}=C_E R^{-5/24},\qquad
     \frac{A}{\kappa R}\le A^{-5/16}=o(1).
     \]
We now compare the matrix $\mathsf C_m$ of lens values
and derivatives with the free fundamental matrix
$\mathsf F_m$. At $r=R$, we have
$p_m(1)^2=1-(m^2-1/4)/(\kappa^2R^2)$, and direct
	differentiation gives
\[
     L_{m,s}'(R)=e^{-siRS_m(2)}
      \left(si\kappa p_m(1)^{1/2}+R^{-1}B_s'(1)\right).
     \]
 We remove the phases from the columns and divide the derivative row by $\kappa$. Both matrices are then
close to the same fixed matrix. We define
	\[
\begin{gathered}
      \mathsf D_m\coloneqq\operatorname{diag}(e^{i\kappa R},e^{-i\kappa R}),\qquad
	\widehat{\mathsf C}_m\coloneqq\mathsf C_m(R)
\operatorname{diag}(e^{iRS_m(2)},e^{-iRS_m(2)}),\\
    \omega_R\coloneqq\frac{1+N^2+A}{\kappa R}=o(1),\qquad
  \mathsf A\coloneqq\begin{pmatrix}1&1\\i&-i\end{pmatrix},\\
E_F\coloneqq\operatorname{diag}(1,\kappa^{-1})\mathsf F_m(R)\mathsf D_m^{-1}
-\mathsf A,\qquad
  E_C\coloneqq\operatorname{diag}(1,\kappa^{-1})\widehat{\mathsf C}_m-\mathsf A.
\end{gathered}
   \]
      The free-wave and lens estimates give
   $\|E_F\|+\|E_C\|\le C\omega_R$. Since $\mathsf A$ is
   a fixed invertible $2\times2$ matrix, $(\mathsf
    A+E_F)^{-1}$ is bounded uniformly for large $b$.
	Therefore
 \[
    \begin{aligned}
      \mathsf D_m\mathsf F_m(R)^{-1}\widehat{\mathsf C}_m
	&=(\mathsf A+E_F)^{-1}(\mathsf A+E_C)\\
    &=I+(\mathsf A+E_F)^{-1}(E_C-E_F)=I+O(\omega_R).
     \end{aligned}
\]
  We express the lens solution in the free basis by defining
   \[
     \binom{d_m}{\overline{d_m}}
  \coloneqq\mathsf F_m(R)^{-1}\mathsf C_m(R)
\binom{\gamma_m}{\overline{\gamma_m}}.
\]
 The two coefficients are conjugates because both pairs
of scalar waves have conjugate values and derivatives.
 The first row of the preceding matrix estimate gives
	\[
   e^{i\kappa R}d_m=e^{-iRS_m(2)}\gamma_m+O(\omega_R)|\gamma_m|.
\]
Taking absolute values, we obtain
      \[
(1-C\omega_R)|\gamma_m|\le|d_m|
      \le(1+C\omega_R)|\gamma_m|,\qquad
     |d_m|\asymp|\gamma_m|\ge c.
\]
      We divide by $d_m$ so that the coefficient of the outgoing free wave at $r=R$
    is one. The resulting reflected solution is
\[
\Psi_m(r,\vartheta)\coloneqq\begin{cases}
     d_m^{-1}\bigl(\gamma_me^{-iRS_m(2)}U_{m,+}
+\overline{\gamma_m}e^{iRS_m(2)}U_{m,-}\bigr),&R\le r\le2R,\\
	 d_m^{-1}w_m(r-2R)\mathcal E_mg(2),&r\ge2R.
   \end{cases}
\]
At $r=R$, the value and radial derivative of $\Psi_m$
 agree with those of the following free solution in the
transformed variables:
      \begin{equation*}
   (F_m^++\rho_mF_m^-)e_m,\qquad \rho_m\coloneqq\overline{d_m}/d_m,\qquad
	|\rho_m|=1.
    \end{equation*}
	For its residual $r_m$, the estimates of
   Lemma~\ref{lem:lens-profiles} and \S\ref{sec:barrier-reflection}
give $\|r_m\|\le\delta_P$. We also need a lower bound for the $L^2$
    mass of $\Psi_m$. The approximation to the ground
state gives
	\[
  q_s=p_m^{-1/2}g+O_P(\epsilon\log(2+A)),\qquad
\|q_s'\|\le C_PA,\qquad |\gamma_m/d_m|\asymp1.
  \]
The two diagonal terms in the squared norm therefore give $2|\gamma_m/d_m|^2R\int_1^2p_m(t)^{-1}\,dt+O_P(R\epsilon\log(2+A))$. The
  remaining term is oscillatory. We set $a_m(t)\coloneqq\langle q_-(t),q_+(t)\rangle$ and integrate by parts. Since $|a_m|\le C$ and
$|a_m'|+|k_m'/k_m|\le CA$,
  \[
R\int_1^2e^{2iRS_m}a_m\,dt
=\left(\frac{e^{2iRS_m}a_m}{2ik_m}\right)_1^2
   -\int_1^2e^{2iRS_m}\left(\frac{a_m}{2ik_m}\right)'dt
     =O_P(A/\kappa).
\]
    Consequently,
 \begin{equation*}
  \begin{aligned}
	 \int_R^{2R}\|\Psi_m(r,\cdot)\|^2\,dr
 &=2|\gamma_m/d_m|^2R\int_1^2p_m(t)^{-1}\,dt
      +O_P(R\epsilon\log(2+A)+A/\kappa),\\
	 c_E R&\le\int_R^{2R}\|\Psi_m(r,\cdot)\|^2\,dr\le C_E R.
     \end{aligned}
\end{equation*}
  Here and below, norms of transformed functions use the measure $dr\,d\vartheta$. The constants are
      uniform for $|m|\le N$, since $A/(\kappa R)=o(1)$. The scalar reflector gives a finite $H^2$ norm
for the part beyond $2R$; its size is not needed in the mass estimate. Radial differentiation and multiplication by the angular potentials preserve the Fourier modes
$m+M\Z$. Thus the construction also gives the asserted angular Fourier support.
   \end{proof}

  \subsection{Matching the solutions and estimating the error}\label{sec:planar-matching}
We now complete the proof of Proposition~\ref{prop:planar-block} using the reflected waves from Lemma~\ref{lem:reflected-waves}. To
solve for general $f\in L^2(B_b)$, we must control the coefficients that match these waves to the solution generated by $f$. We
	obtain a bound uniform in the angular mode by excluding a small neighborhood of the spectrum of a single Dirichlet operator.

\begin{proof}[Proof of Proposition~\ref{prop:planar-block}]
 We let $H_b^D$ be the realization of $-\Delta+V_b$ on $B_{3R}$ with Dirichlet domain $H^2(B_{3R})\cap H_0^1(B_{3R})$, and define
\[
   Z_b\coloneqq\{\mu\in[R^{-1/4},E]:
   \dist(\mu,\sigma(H_b^D))<R^{-L}\}.
  \]
Since $V_b\ge0$, the min--max principle gives, with all angular multiplicities included,
   \[
\begin{gathered}
    \rank\1_{(-\infty,E+1]}(H_b^D)\le CR^2,\qquad
 |Z_b|\le CR^{2-L}=CR^{-T-1},\\
\|(H_b^D-\mu)^{-1}\|\le R^L
\quad\text{for }\mu\in[R^{-1/4},E]\setminus Z_b.
  \end{gathered}
	\]
We exclude these energies once for all $f\in L^2(B_b)$ and all angular modes $m$.

 Recall the reflection phase $\rho_m$ and reflected wave $\Psi_m$ from Lemma~\ref{lem:reflected-waves}. We set $j_m\coloneqq e^{-2i\alpha_m}$ and
      $h_m\coloneqq1-j_m\rho_m$. The number $h_m$ measures whether the reflected wave can be continued regularly to the origin. If it were too small, the continuation
below would be an accurate Dirichlet quasimode at an energy outside $Z_b$. We choose a fixed smooth $\chi$, zero on $(-\infty,1/2]$ and one on $[3/4,\infty)$, and set
    $\chi_R(r)\coloneqq\chi(r/R)$. In the transformed polar coordinates define
      \[
    \widetilde v_m(r,\vartheta)\coloneqq
	 \begin{cases}
(\rho_m(F_m^-+j_mF_m^+)+h_m\chi_RF_m^+)e_m(\vartheta),&r\le R,\\
	 \Psi_m(r,\vartheta),&r\ge R.
    \end{cases}
    \]
The values and radial derivatives match at $r=R$, because $\chi_R=1$ there and
 \[
	\rho_m(F_m^-+j_mF_m^+)+h_mF_m^+=F_m^++\rho_mF_m^-.
 \]
      At $r=2R$, the values and radial derivatives match by the definition of $\mathsf C_m(2R)$ and $\gamma_m$. At the origin, the identity
\[
	r^{-1/2}(F_m^-+j_mF_m^+)e_m
 =\sqrt{2\pi\kappa}\,e^{-i\alpha_m}
J_{|m|}(\kappa r)e_m
\]
 gives a smooth Cartesian function. The term with $\chi_R$ vanishes near the origin. Thus the physical function $v_m\coloneqq
 r^{-1/2}\widetilde v_m$ belongs to $\mathcal D(H_b^D)$ and vanishes outside $B_{9R/4}$. Since the values and first derivatives
	match, the distributional Laplacian has no delta terms at the interfaces.

	On $R/2\le r\le R$, the uniform Hankel bounds give
 \[
  |F_m^+(r)|\le C,\qquad |(F_m^+)'(r)|\le C\kappa,
	 \qquad \frac{N^2+1}{\kappa R}=o(1).
 \]
    Since $V_b=0$ there, the only residual in the inner region is the cutoff commutator:
\[
 \begin{aligned}
    \|[-\partial_r^2,\chi_R]F_m^+\|_{L^2(dr)}
 &\le C(R^{-2}R^{1/2}+R^{-1}\kappa R^{1/2})\\
&\le C\kappa R^{-1/2},\qquad \kappa R\ge R^{7/8},\\
  \|(H_b^D-\mu)v_m\|
      &\le\delta_P+C|h_m|\kappa R^{-1/2},\qquad
     \|v_m\|\ge cR^{1/2}.
 \end{aligned}
      \]
   All constants are independent of $m$ and $N$. A small $h_m$ would make $v_m$ a Dirichlet quasimode. For $\mu\notin Z_b$,
\[
R^{-L}\|v_m\|\le\|(H_b^D-\mu)v_m\|,
	\]
    and therefore
	\[
cR^{1/2-L}\le\delta_P+C|h_m|\kappa R^{-1/2}.
 \]
  The order choice \eqref{eq:planar-expansion-order} ensures
    \[
\begin{aligned}
R^{L-1/2}\delta_P
&\le C_{P,E}A^{r_*(L-1)+1/2-\ell P/2}
+R^{L-1/2}\varepsilon_{{\rm bar},E}(A)\longrightarrow0.
     \end{aligned}
\]
   Absorbing this term gives, simultaneously for every
$|m|\le N$,
      \begin{equation*}
|1-j_m\rho_m|\ge c\kappa^{-1}R^{1-L},\qquad
    \max_{|m|\le N}|1-j_m\rho_m|^{-1}\le C\kappa R^{L-1}.
   \end{equation*}

	For $f\in L^2(B_b)$, we first solve the free equation using
 \[
w_f(x)\coloneqq\frac i4\int_{B_b}H_0^{(1)}(\kappa|x-y|)f(y)\,dy.
\]
 This solves $(-\Delta-\mu)w_f=\iota_bf$ with outgoing boundary conditions. Let $\mathbf b_f\coloneqq((b_f)_m)_{m\in\Z}$ be its coefficients: for $r>b$,
\[
\sqrt r\,w_f=\sum_{m\in\mathbb Z}(b_f)_mF_m^+(r)e_m(\vartheta),
  \qquad \|\mathbf b_f\|_{\ell^2}\le Cb\kappa^{-1/2}\|f\|.
  \]
For $\omega\coloneqq(\cos\vartheta,\sin\vartheta)$ set
    $a_f(\vartheta)\coloneqq\int_{B_b}e^{-i\kappa\omega\cdot
 y}f(y)\,dy$ and $c_\kappa\coloneqq
      ie^{-i\pi/4}/\sqrt{8\pi\kappa}$. For fixed
	$b,\kappa>0$,
    \[
\begin{gathered}
      e^{-i\kappa r}\sqrt r\,w_f(r,\cdot)\xrightarrow[r\to\infty]{L^2(\mathbb S^1)}
c_\kappa a_f,\\
      \sum_m(b_f)_me_m=c_\kappa a_f,\qquad \|a_f\|_2\le Cb\|f\|.
 \end{gathered}
  \]
  Taking each fixed Fourier coefficient in this limit identifies $(b_f)_m$ with the coefficient of $c_\kappa a_f$. Orthonormality of $(e_m)$ then gives the estimate explicitly:
\[
  \sum_m|(b_f)_m|^2=|c_\kappa|^2\|a_f\|_2^2
	 \le Cb^2\kappa^{-1}\|f\|_2^2.
    \]
	We add a regular free solution to match the reflected wave at $R$. The outgoing coefficient is then $(b_f)_m+j_m\xi_m$, while the incoming coefficient is $\xi_m$. Matching them to $\tau_m(1,\rho_m)$ requires $\tau_m=(b_f)_m+j_m\xi_m$ and $\xi_m=\rho_m\tau_m$. We therefore define the vectors
$\boldsymbol\tau,\boldsymbol\xi\in\C^n$ by
 \[
\begin{gathered}
	\tau_m\coloneqq\frac{(b_f)_m}{1-j_m\rho_m},\quad
\xi_m\coloneqq\rho_m\tau_m,\quad \tau_m=(b_f)_m+j_m\xi_m\quad\text{for }|m|\le N,\\
   \|\boldsymbol\tau\|+\|\boldsymbol\xi\|
      \le Cb\kappa^{1/2}R^{L-1}\|f\|.
     \end{gathered}
	 \]
The corresponding regular free correction is
  \[
   \sqrt r\,W_{\boldsymbol\xi}
\coloneqq\sum_{|m|\le N}\xi_m(F_m^-+j_mF_m^+)e_m.
	\]
The identity $H_j^{(1)}+H_j^{(2)}=2J_j$ makes the
  correction smooth at the origin.

	We write $P_{>N}$ for the projection onto the Fourier modes $|m|>N$. We choose a fixed smooth cutoff $\zeta_{\rm cut}$, one on $(-\infty,1/4]$ and zero on $[1/2,\infty)$, and take $\zeta_R(r)\coloneqq\zeta_{\rm cut}(r/R)$. We
   have matched the modes $|m|\le N$. We cut off the remaining modes inside $R/2$, before the lens begins. We define
  \[
    u_{\rm in}\coloneqq w_f+W_{\boldsymbol\xi}-(1-\zeta_R)P_{>N}w_f.
  \]
   For $r\in[R/4,R/2]$ and $|y|\le b$, fix a strip
$|\Im z|\le v_0$ with $v_0>0$ independent of $b$. We
  define the complex distance by
    \[
   \widetilde r(z)\coloneqq\sqrt{r^2-2r(\cos z,\sin z)\cdot y+|y|^2},
\]
	using the square root near $+r$. Since the expression under the root is $r^2+O(rb)$ and $b/r\to0$, this branch is analytic and periodic in the fixed strip. It satisfies
  \[
	\widetilde r=r+O(b),\qquad \Re\widetilde r\ge r/2,\qquad
|\Im\widetilde r|\le Cb,\qquad \partial_r\widetilde r=1+O(b/r).
	\]
	The kernel and its radial derivative are therefore bounded by $C_E\kappa^{-1/2}r^{-1/2}e^{C_Eb}$, since $\kappa r\ge cR^{7/8}$. For a periodic analytic function $a$ bounded by $M_a$ on this strip,
	\[
	|\widehat a(m)|\le C M_ae^{-v_0|m|},\qquad
 \sum_{|m|>N}|\widehat a(m)|^2\le C M_a^2e^{-2v_0N}.
	\]
We apply this Fourier estimate to
$\partial_r^jH_0^{(1)}(\kappa\widetilde r(z))$ and integrate against
 $f$. Cauchy--Schwarz on $B_b$ then gives, for $j=0,1$,
  \[
      \|P_{>N}\partial_r^jw_f(r,\cdot)\|_{L^2(\mathbb S^1)}
  \le Cb\kappa^{-1/2}r^{-1/2}e^{C_Eb-cN}\|f\|.
\]
This cutoff changes the equation only where $\zeta_R$ varies. In the original polar coordinates, the resulting error is
\[
\begin{aligned}
    [-\Delta,\zeta_R]P_{>N}w_f
  &=-2\zeta_R'P_{>N}\partial_rw_f
   -(\zeta_R''+r^{-1}\zeta_R')P_{>N}w_f,\\
\|[-\Delta,\zeta_R]P_{>N}w_f\|
	 &\le Cb\kappa^{-1/2}e^{C_Eb-cN}R^{-1/2}\|f\|.
\end{aligned}
	 \]
	Fixing constants in $N=\lceil C_0b+C_1\log R\rceil$ and choosing $cC_0>C_E$ and $cC_1>B+1$ give
    \[
  \begin{gathered}
 e^{C_Eb-cN}\le e^{-(cC_0-C_E)b}R^{-cC_1},\qquad
  b\kappa^{-1/2}R^{-1/2}\le CR^{1/3+1/16-1/2}=CR^{-5/48},\\
Cb\kappa^{-1/2}e^{C_Eb-cN}R^{-1/2}\le CR^{-cC_1-5/48}=o(R^{-B}).
 \end{gathered}
\]
On $R/2\le r\le R$, the transformed inner solution satisfies $r^{1/2}u_{\rm in}=\sum_{|m|\le N}\tau_m(F_m^++\rho_mF_m^-)e_m$. The value and radial derivative at $R$ match those of $\sum_{|m|\le N}\tau_m\Psi_m$. We therefore define
\[
 Q_b(\mu)f(r,\vartheta)\coloneqq
\begin{cases}u_{\rm in}(r,\vartheta),&r\le R,\\
r^{-1/2}\sum_{|m|\le N}\tau_m\Psi_m(r,\vartheta),&r\ge R.
\end{cases}
  \]
    This function is supported in $B_{9R/4}$. At every joining radius,
    \[
[Q_bf]=[\partial_rQ_bf]=0
      \ \Longrightarrow\ \Delta Q_bf=(\Delta Q_bf)_{\rm pieces}\in L^2
\ \Longrightarrow\ Q_bf\in H^2(\R^2).
\]
  The bounds on each piece follow from local elliptic
   regularity and the estimates on $H^2(\mathbb S^1)$
proved in \S\ref{sec:lens-solutions}. To return from
   the periodic angular variable to $\vartheta$, we use
\[
D_\vartheta\mathcal E_mz=\mathcal E_m(MD+m)z,\qquad
    \|\mathcal E_mz\|_{H^2}\le C(1+M^2+N^2)\|z\|_{H^2}\quad\text{for }|m|\le N.
    \]
    Thus $Q_b(\mu):L^2(B_b)\to H^2(\R^2)$ is bounded for each fixed $b,\mu$. Its $H^2$ norm may grow with $M,N$ and with $e^{C_Ed_A}$. We need the bound for the error in
the equation to remain uniform, and we estimate that error next.

  Each $r_m$ has Fourier support in $m+M\Z$. Radial differentiation and multiplication by $W_a(M\vartheta)$ or $F_A(M\vartheta)$
      preserve this set. Since $M>2N$, these sets are disjoint for $|m|\le N$. The orthogonality identity
    \eqref{eq:planar-ladder-orthogonality} therefore
	gives
	\[
\left\|\sum_{|m|\le N}\tau_m{r_m}\right\|^2
     =\sum_{|m|\le N}|\tau_m|^2\|{r_m}\|^2
   \le\delta_P^2\|\boldsymbol\tau\|^2.
\]
We have estimated the sum by the Euclidean norm of $\boldsymbol\tau$, with a constant independent of $N$. Combining this with the coefficient bound above, the polynomial
part of the error is at most $\|f\|$ times
\[
      \begin{aligned}
     C_{P,E}b\kappa^{3/2}R^{L-3/2}A^{1/2}[\epsilon\log(2+A)]^P
 &\le C_{P,E}A^{r_*L-5/4-\ell P/2}
 =o(R^{-B}),
 \end{aligned}
\]
where $b\le CA^{1/2}$, $\kappa\le\sqrt E$, and
 \eqref{eq:planar-expansion-order} were used. The barrier
contribution is bounded by
\[
    C_EA^{r_*L-1}\varepsilon_{{\rm bar},E}(A)\|f\|
   =o(R^{-B})\|f\|.
 \]
Together with the bound on the higher Fourier modes, this proves Proposition~\ref{prop:planar-block} for sufficiently large $b$.
    \end{proof}

\subsection{The planar operator}\label{sec:planar-operator}

	We now sum the annular potentials and apply
      Proposition~\ref{prop:assembly}. We first choose the energy
ceilings so that one sum works on every compact interval of positive
  energies. This choice will also be used for the axial barriers in
	\S\ref{sec:dimensions}.

For every integer $k\ge1$, construct the planar
barriers with ceiling $E=k$. Choose a strictly
increasing sequence of thresholds $T_k\to\infty$ with
   the following properties. For every $b\ge T_k$, the
potential $V_{b,k}$ has the height, support, and cost
	bounds stated above, and
Proposition~\ref{prop:planar-block} holds for all
   positive integer orders $B,T\le k+650$, with error
	at most $R^{-B}$ and excluded measure at most
  $R^{-T}$. There are only finitely many requirements at each $k$. To obtain the
 unit coefficient in the measure estimate, apply the proposition with one larger
      order and increase $T_k$ to absorb its fixed constant.

 We include the axial estimates of Lemma~\ref{lem:axial-inverse} in the same choice of thresholds. This lemma is independent of the planar background $U$. Let $C_k,c_k>0$ be its constants at ceiling $k$ and exceptional-set order $L=8$, with $C_k$ enlarged to include the localized derivative estimates used in
the cylinder construction. Increase $T_k$ so that this lemma applies at $Z=b^{40}$, with its stated support bounds, for every $b\ge T_k$. Also require
\begin{equation}\label{eq:diagonal-axial-bounds}
  C_k\le b,\qquad
 C_k Z^{12}\exp\bigl(-c_k(\log Z)^2\bigr)\le Z^{-152}.
    \end{equation}
These requirements can be met for each $k$, since $c_k>0$ and the exponential decays faster than every power of $Z$. No lower bound
 on $c_k$ uniform in $k$ is needed. Define the nondecreasing function
\begin{equation}\label{eq:planar-ceiling}
K(b)\coloneqq\max\{k\ge1:T_k\le b\},\qquad b\ge T_1.
\end{equation}
Then $K(b)\to\infty$. At each source radius $b$ use the
planar potential $V_b\coloneqq V_{b,K(b)}$; in the
 cylinder construction, use the axial barrier with the same ceiling $K(b)$. At any fixed energy ceiling and for any fixed error orders, the required estimates hold at every sufficiently large stage. The radii retain the
 recursion below; the ceiling increases only when that recursion
reaches its next threshold.

We write $\rho\in\R^2$ for the planar coordinate. Choose $\beta_2+1>15$ and $Q_2>27$. Starting from a large $b_1\ge T_1$, let $R_n$ be the radius corresponding to $b_n$ in \eqref{eq:planar-parameters}, and set
    $b_{n+1}\coloneqq3R_n$. Write 
   \[
U\coloneqq\sum_nV_{b_n},\qquad h_0\coloneqq-\Delta_\rho+U
 \quad\text{on }H^2(\R^2),\qquad W_2\coloneqq\langle\rho\rangle^{-2Q_2}.
      \]
The supports are disjoint and the sum is locally
	finite. Moreover, $10R_n\le R_{n+1}\le CR_n^3$ once $b_1$ is large.

   \begin{proposition}\label{prop:planar-quantitative}
  The function $U$ is smooth and $0\le U\le4$. Taking $b_1$ sufficiently large makes $\mathcal
I_{2,\sigma}(U)$ arbitrarily small. For every fixed $E,B,T>0$, all sufficiently large $n$ have the
  conclusion of Proposition~\ref{prop:planar-block} on $[R_n^{-1/4},E]$ for the actual summand $V_{b_n}$,
	with error at most $R_n^{-B}$ and excluded measure at most $R_n^{-T}$. For almost every $s>0$,
      $h_s\coloneqq h_0+sW_2$ has pure point spectral type on $(0,\infty)$ with finite eigenvalue
multiplicities, zero kernel, and $\sigma(h_s)=\sigma_{\mathrm{ess}}(h_s)=[0,\infty)$. Taking $b_1$ large
    and then choosing a sufficiently small $s>0$ for which this conclusion holds makes both weighted integrals of $V_2\coloneqq U+sW_2$ arbitrarily small.
      \end{proposition}

 \begin{proof}
The quantitative statement follows from the threshold construction: eventually $K(b_n)\ge E$ and $K(b_n)+650$ exceeds both requested orders, rounded upwards to integers.
The block chosen at that stage therefore covers $[R_n^{-1/4},E]$ with the asserted estimates.

   On each compact $J\subset(0,\infty)$, apply this
 statement with $B\coloneqq\beta_2$ and $T\coloneqq1$
for all sufficiently large $n$, when also $R_n^{-1/4}<\inf J$. We take $Z_n\coloneqq J$ at every earlier stage. Proposition~\ref{prop:assembly} then applies with
\[
d=2,\qquad p\coloneqq3,\qquad \alpha\coloneqq3,\qquad
     \beta\coloneqq\beta_2,\qquad Q\coloneqq Q_2,
 \]
by our choices of $\beta_2$ and $Q_2$. The countable exhaustion $[1/j,j]$, $j\ge1$, gives a set of full measure of values of $s$ for
which the conclusion holds throughout $(0,\infty)$, for the same $U$. The cost constants in \eqref{eq:planar-cost} are independent
of the ceiling. Hence those weighted integrals are summable and their sum tends to zero as $b_1\to\infty$. Since all annular
 supports lie outside $B_1$, $\mathcal I_{2,0}(U)\le\mathcal I_{2,\sigma}(U)$. Both weighted integrals of $W_2$ are finite, so
choosing $s>0$ sufficiently small outside the exceptional set keeps both integrals arbitrarily small.

If $h_su=0$, then $\|\nabla u\|^2+\langle(U+sW_2)u,u\rangle=0$, so $u=0$.
      Nonnegativity and the positive-energy conclusion thus give $P_c(h_s)=0$. The balls $B_{R_n}(5R_n\mathbf e_1)$ avoid the supports of all annular potentials, and
    $W_2$ tends uniformly to zero there. Their cutoff plane waves give Weyl sequences at every nonnegative energy.
  \end{proof}

\section{Cylindrical confinement and higher dimensions}\label{sec:dimensions}

 We use the annular potential $U$ from
Proposition~\ref{prop:planar-quantitative} as the
 background, without the term $sW_2$.

      \subsection{A compact inverse at low energy}

 Here $h_0=-\Delta_\rho+U$ is the planar operator from \S\ref{sec:planar-operator}. When we separate the
    remaining coordinate, we subtract a nonnegative axial eigenvalue from the energy. The resulting
 planar energy may be close to zero or arbitrarily negative. We therefore need an approximate inverse
that remains accurate in this range. At these energies, the angular lens alone makes the solution small
near the outer boundary.

\begin{lemma}\label{lem:low-inverse}
	There are constants $c,C>0$ such that, for every sufficiently large planar radius $S=R_m$ and every real $\mu\le cS^{-2}$, there is
a bounded map
  \[
Q_S^{\mathrm{in}}(\mu):L^2(B_{S/2})\longrightarrow H^2(\R^2),
	 \qquad \supp Q_S^{\mathrm{in}}(\mu)f\subset B_{2S},
\]
 with
  \begin{equation}\label{eq:low-inverse}
   \|(h_0-\mu)Q_S^{\mathrm{in}}(\mu)-\iota\|
      \le Ce^{-cS^{1/6}}.
\end{equation}
    The constant is uniform for all these energies.
   \end{lemma}

  \begin{proof}
  We use the parameters $A,M$ at radius $S$. Where the lens parameter equals $A$,
\[
     D_\vartheta^2+M^2W_A(M\vartheta)\ge cM^2A.
\]
 The diamagnetic inequality \cite[\S9]{SimonDia} holds for every real $\zeta$, so
\[
	 \inf\sigma(K_{A,\zeta})\ge cA\quad\text{for }\zeta\in\R,\qquad
	 D_\vartheta^2+M^2W_A(M\vartheta)
\simeq\bigoplus_{m=0}^{M-1}M^2K_{A,m/M}\ge cM^2A.
\]
The lower bound follows from \eqref{eq:planar-angular-lower}. The schedule equals
$A$ on a terminal interval of length comparable to $A^{-1}$ in the
variable $r/S$. We therefore choose a fixed $\epsilon_0>0$ small
enough that
\[
      \mathcal C_S\coloneqq
 \{(2-\epsilon_0/A)S<|\rho|<(2-\epsilon_0/(2A))S\}
\]
lies inside that interval. We set
    \[
  v_S\coloneqq M^2A/S^2\asymp A^{\ell+1-r_*},\qquad r_*=\tfrac32.
    \]
   We let $h_D$ be the Dirichlet restriction of $h_0$
      to $B_{2S}$ and write $q_D$ for its quadratic
form. Since $U\ge0$, $h_D\ge g_S I$, where
$g_S\coloneqq c_0S^{-2}$ for a fixed $c_0>0$. We
 integrate the angular inequality in $r$ to obtain, for
    $w\in H_0^1(B_{2S})$, ${q_D}(w)\ge
  cv_S\|\1_{\mathcal C_S}w\|^2$. We add half of each bound to obtain
      \[
      {q_D}(w)\ge\tfrac12g_S\|w\|^2
+\tfrac12cv_S\|\1_{\mathcal C_S}w\|^2.
  \]
We choose a radial Lipschitz weight $\Phi$, equal to zero for $|\rho|\le(2-\epsilon_0/A)S$, increasing across the annulus $\mathcal C_S$, and constant for $|\rho|\ge(2-\epsilon_0/(2A))S$. We require $|\nabla\Phi|^2\le
    cv_S/4$. We take $\Phi'=c_1\sqrt{v_S}$ on a fixed middle portion of $\mathcal C_S$ and zero elsewhere, with $c_1>0$ small enough. Then
\[
\Phi=0\ \text{before }\mathcal C_S,\qquad
\Phi\ge c(S/A)\sqrt{v_S}\ \text{after }\mathcal C_S,\qquad
|\nabla\Phi|^2\le cv_S/4.
	\]
   For $\supp f\subset B_{S/2}$ and $\mu\le g_S/4$, we
	set $u\coloneqq(h_D-\mu)^{-1}f$ and $w\coloneqq e^\Phi
	u$. We apply the Agmon identity \cite[Proposition~6.1]{TZ} and the preceding form bound to obtain
	\[
	\tfrac14g_S\|w\|^2+\tfrac14cv_S\|\1_{\mathcal C_S}w\|^2
     \le{q_D}(w)-\mu\|w\|^2-\|w\nabla\Phi\|^2
=\Re\langle w,f\rangle.
    \]
We apply Cauchy--Schwarz to obtain
\[
 \|w\|\le CS^2\|f\|,\qquad
v_S\|\1_{\mathcal C_S}w\|^2+\|w\nabla\Phi\|^2\le CS^2\|f\|^2.
	 \]
    Since $U\ge0$ and $\mu\le g_S/4$, the same identity gives
 \[
      \|\nabla w\|^2\le{q_D}(w)
=\Re\langle w,f\rangle+\mu\|w\|^2+\|w\nabla\Phi\|^2
 \le CS^2\|f\|^2.
\]
Using $e^\Phi\nabla u=\nabla w-w\nabla\Phi$, we obtain
 \[
	S^{-2}\|e^\Phi u\|+S^{-1}\|e^\Phi\nabla u\|\le C\|f\|.
\]
The constants are uniform for $\mu\le g_S/4$. The lower bound only improves when
$\mu$ decreases below zero.

We choose a radial cutoff $\theta_S$, equal to one on $B_{S/2}$ and zero near
	$\partial B_{2S}$, with derivatives supported where $\Phi\ge c(S/A)\sqrt{v_S}$ and
    bounded by $C_j(A/S)^j$ for $j\le2$. We define
 $Q_S^{\mathrm{in}}(\mu)f\coloneqq\theta_Su$ and extend it by zero. Then
   \[
\begin{gathered}
(h_0-\mu)Q_S^{\mathrm{in}}(\mu)f-f=[-\Delta,\theta_S]u,\\
      \|[-\Delta,\theta_S]u\|
\le C(A/S)\|\nabla u\|_{\supp\nabla\theta_S}
  +C(A/S)^2\|u\|_{\supp\Delta\theta_S} \le C(A+A^2)e^{-c(S/A)\sqrt{v_S}}\|f\|.
   \end{gathered}
\]
The exponent satisfies
\[
(S/A)\sqrt{v_S}\asymp A^{(r_*+\ell-1)/2}
      =S^{1/6+\ell/3}.
\]
  Absorbing the polynomial factor proves \eqref{eq:low-inverse}.
  Interior elliptic regularity gives $Q_S^{\mathrm{in}}(\mu)f\in H^2$. The
 constants use only the angular lens and nonnegativity of the other blocks, and are independent of their
 energy ceilings.
\end{proof}

   \subsection{A Dirichlet cylinder}

      We write $x=(\rho,z)\in\R^2\times\R$ and set $H_{\mathrm{cyl}}\coloneqq-\Delta_\rho-\partial_z^2+U(\rho)$. We use the annular radii $R_k$ from \S\ref{sec:planar-operator}. For $f$ supported in $B_b$, we choose $n\coloneqq\min\{k:R_k\ge b^{10}\}$ and set
\begin{equation}\label{eq:cylinder-scales}
      S\coloneqq R_n,\qquad b^{10}\le S\le Cb^{30},\qquad
    Z\coloneqq b^{40},\qquad w_b\coloneqq(\log Z)^2.
  \end{equation}
      The axial barriers consist of finitely many periodic pieces. We
first give the approximate inverse used near the two ends of the
	cylinder.

 \begin{lemma}\label{lem:axial-inverse}
Fix $E>0$ and $L>2$. For all sufficiently large $Z$,
set $w_Z=(\log Z)^2$ and $I_Z=(-Z-w_Z,Z+w_Z)$. There is
      an even smooth potential $v_Z$, with
  \[
   0\le v_Z\le3,\qquad
 \supp v_Z\subset\{z:Z<|z|<Z+w_Z\},
\]
      and a measurable set $\mathcal E_Z\subset[0,E]$
	of measure
\begin{equation}\label{eq:axial-exceptional-set}
      |\mathcal E_Z|\le C_E Z^{1-L},
\end{equation}
    such that for every real $\nu\le E$ outside $\mathcal E_Z$ there is a bounded linear map
\[
      Q_Z^{\mathrm{ax}}(\nu):L^2((-Z,Z))\longrightarrow H^2(\R)
      \]
with
	\begin{equation}\label{eq:axial-local-inverse}
\begin{gathered}
   \supp Q_Z^{\mathrm{ax}}(\nu)f
  \subset[-Z-w_Z/3,Z+w_Z/3],\qquad
    \|Q_Z^{\mathrm{ax}}(\nu)f\|\le C_E Z^L\|f\|,\\
   \|(-\partial_z^2+v_Z-\nu)Q_Z^{\mathrm{ax}}(\nu)f-\iota f\|
	\le C_E Z^L e^{-c_Ew_Z}\|f\|.
    \end{gathered}
\end{equation}
Here $\iota$ denotes extension by zero. The constants are uniform also as $\nu\to-\infty$. In particular, all these approximate inverses vanish on the same outer collar
	\[
      \mathcal T_Z=\{z:Z+w_Z/2<|z|<Z+3w_Z/4\}\subset I_Z.
	\]
 The potential $v_Z$ does not depend on $L$.
      \end{lemma}

	\begin{proof}
      We use Lemma~\ref{lem:periodic-gaps} and first describe the positive half of the potential. Let $m=n_E+1$. For each $i$, choose a positive integer $n_i$ so that
 \[
L_i=n_iT_i\asymp_E w_Z,\qquad
 \sum_{i=0}^{n_E} L_i\le w_Z/8.
\]
	For example, $n_i=\lfloor w_Z/(8mT_i)\rfloor$ works for large $Z$. The $i$th periodic part starts at
\[
a_i=Z+3+\sum_{h<i}(L_h+2)
  \]
	and equals $q_i(z-a_i)$ on $[a_i,a_i+L_i]$. Connect successive periodic parts by
smooth convex interpolation on intervals of length two, and connect the first
constant part to zero inside $(Z,Z+3)$. The endpoint
\[
D_i=a_i+L_i
\]
 of each periodic part has the prescribed phase zero modulo $T_i$. All these endpoints lie below
$Z+w_Z/4$ for large $Z$. Continue the last periodic function beyond $D_{n_E}$ to $Z+3w_Z/4$ and smoothly
	cut it off before $Z+7w_Z/8$. Reflection across zero gives the required even
smooth $v_Z$, bounded between zero and three.

   Let $h_{Z,i}$ be the Dirichlet realization of $-\partial_z^2+v_Z$ on $(-D_i,D_i)$, and write
\[
  \mathcal E_Z=\bigcup_{i=0}^{n_E}
	 \{\nu\in[0,E]:\dist(\nu,\sigma(h_{Z,i}))<Z^{-L}\}.
\]
  Nonnegativity and comparison with the Dirichlet Laplacian give
   \[
  \rank\1_{(-\infty,E+1]}(h_{Z,i})\le C_EZ.
 \]
There are finitely many $h_{Z,i}$, which proves \eqref{eq:axial-exceptional-set}.
   For $\nu\in[0,E]\setminus\mathcal E_Z$,
\begin{equation}\label{eq:axial-dirichlet-inverse}
  \|(h_{Z,i}-\nu)^{-1}\|\le Z^L.
  \end{equation}
Also $h_{Z,i}\ge cZ^{-2}$, so for all $\nu<0$,
\begin{equation}\label{eq:negative-axial-inverse}
\|(h_{Z,i}-\nu)^{-1}\|\le CZ^2\le C Z^L.
	 \end{equation}
The constants in the last bound are independent of the negative energy.

We first suppose $\nu\le\min\{2,E\}$, and use $i=0$. Extend $f\in
	L^2((-Z,Z))$ by zero and set $u=(h_{Z,0}-\nu)^{-1}f$. On the two end
intervals $a_0<|z|<D_0$, the potential equals three. Choose a fixed
$a>0$ with $a^2<1/2$, and define on $(-D_0,D_0)$
   \[
 \Phi(z)=a(|z|-a_0)_+,\qquad W=e^\Phi u,
	 \qquad G=\{z:a_0<|z|<D_0\}.
\]
   Testing the weak equation against $e^{2\Phi}u$ gives
  \[
      \int|W'|^2+\int(v_Z-\nu-|\Phi'|^2)|W|^2
=\Re\langle u,f\rangle.
  \]
On the complement of $G$, $\Phi=0$ and $v_Z\ge0$; on $G$, $v_Z-\nu-a^2\ge1-a^2>1/2$. Consequently
   \[
 \|W'\|^2+\tfrac12\|\1_GW\|^2
      \le 2\|u\|^2+\|u\|\|f\|
   \le C_E Z^{2L}\|f\|^2.
\]
This estimate only uses $\nu\le2$ and \eqref{eq:axial-dirichlet-inverse}--\eqref{eq:negative-axial-inverse}, so it is uniform down to $\nu=-\infty$. On the middle third
      of each end interval, $\Phi\ge aL_0/3\ge c_Ew_Z$. The identity $e^\Phi u'=W'-\Phi'W$ therefore gives
   \begin{equation}\label{eq:low-axial-collar}
\|u\|_{\mathrm{middle}}+\|u'\|_{\mathrm{middle}}
\le C_E Z^L e^{-c_Ew_Z}\|f\|.
	\end{equation}

We next suppose $2<\nu\le E$ and $\nu\notin\mathcal E_Z$. Choose $i\ge1$ with $\nu\in J_i$ and put
	$u=(h_{Z,i}-\nu)^{-1}f$. On the right end interval, use the coordinate $s=z-a_i\in[0,L_i]$. Let
$v_s,v_u$ be unit stable and unstable eigenvectors of the monodromy, with multipliers $\lambda_s$ and
    $\lambda_u$. These satisfy
    \[
|\lambda_s|<1<|\lambda_u|,\qquad
  \lambda_s\lambda_u=1,\qquad |(v_u)_1|\ge c_E.
  \]
    Write the Cauchy data at $s=0$ as $Av_s+Bv_u$. Since $L_i=n_iT_i$ and $u(L_i)=0$, the endpoint condition gives
    \[
	 A\lambda_s^{n_i}(v_s)_1
    +B\lambda_u^{n_i}(v_u)_1=0,
\qquad |B|\le C_E e^{-2\gamma_E L_i}|A|.
   \]
  The within-period propagators and the monodromy eigenprojections are uniformly bounded on the compact
energy interval. Thus on $L_i/3\le s\le2L_i/3$,
  \[
     |u(a_i+s)|+|u'(a_i+s)|
 \le C_Ee^{-\gamma_E L_i/3}(|u(a_i)|+|u'(a_i)|).
    \]
A local one-dimensional energy estimate followed by the Sobolev
   inequality, applied in a fixed neighborhood of $a_i$, bounds the last Cauchy norm by $C_E\|u\|$. This neighborhood is disjoint from $\supp f$ and lies inside
   $(-D_i,D_i)$ for large $Z$; the equation there is homogeneous with $|v_Z-\nu|\le3+E$. The same estimate holds at the left endpoint, using the outward coordinate $s=-z-a_i$. Integrating the pointwise bounds over the middle
 thirds contributes at most a factor $C_Ew_Z^{1/2}$, which is absorbed into the exponential. By \eqref{eq:axial-dirichlet-inverse},
\begin{equation}\label{eq:high-axial-collar}
    \|u\|_{\mathrm{middle}}+\|u'\|_{\mathrm{middle}}
      \le C_E Z^L e^{-c_Ew_Z}\|f\|.
\end{equation}

   In either case we choose an even smooth cutoff $\chi_i$, equal to one for
	$|z|\le a_i+L_i/3$ and zero for $|z|\ge a_i+2L_i/3$, with its first two
    derivatives supported in the intervening middle thirds and bounded by $C_E$. We
define
    \[
	 Q_Z^{\mathrm{ax}}(\nu)f=\chi_i u
\]
 and extend the result by zero to the line. The cutoff equals one on
  $(-Z,Z)$ and vanishes near the endpoints of the Dirichlet interval.
    Consequently the extension belongs to $H^2(\R)$, and
\[
	(-\partial_z^2+v_Z-\nu)Q_Z^{\mathrm{ax}}(\nu)f-\iota f
	=-2\chi_i'u'-\chi_i''u.
  \]
Equations~\eqref{eq:low-axial-collar} and
\eqref{eq:high-axial-collar} prove the error bound. The
$L^2$ bound follows from the Dirichlet resolvent
	bounds. Finally its support lies strictly below
 $D_i<Z+w_Z/4$ in absolute value, and hence inside the
stated interval and away from the common collar
$\mathcal T_Z$.
 \end{proof}

 We add two disk-shaped barriers near $z=Z$ and $z=-Z$.
	We use the nondecreasing energy ceiling $K(b)$ chosen
	in \S\ref{sec:planar-operator}, and let $v_b$ be the
 potential of Lemma~\ref{lem:axial-inverse} with
      $E=K(b)$ and $Z=b^{40}$. We also choose $A_b\in C_c^\infty(\R^2;[0,1])$, equal to one on $B_{3S}$ and zero outside $B_{4S}$, and define
  \[
V_b^{\mathrm{cap}}(\rho,z)\coloneqq A_b(\rho)v_b(z),\qquad
     \supp V_b^{\mathrm{cap}}
 \subset\{|\rho|\le4S,\ Z<|z|<Z+w_b\}.
 \]
   Since $S/Z\le Cb^{-10}$ and $w_b/Z\to0$, this
support lies in $\{Z<|x|<9Z/4\}$ for large $b$. We will
	solve a Dirichlet problem in a cylinder containing
    $B_b$. The planar inverse controls the solution
  near the curved boundary. The two barriers control it near the ends of the cylinder. We can then cut
    the solution off before it reaches any boundary.

  \begin{proposition}\label{prop:cylinder-block}
Fix $E>0$. For all sufficiently large $b$, there is a measurable set
$\mathcal Z_b\subset[0,E]$, $|\mathcal Z_b|\le C_Eb^{-1}$, such that
      for each $\lambda\in[0,E]\setminus\mathcal Z_b$ there is a
bounded linear map $Q_b(\lambda):L^2(B_b)\to H^2(\R^3)$ satisfying, for every $f\in L^2(B_b)$,
\begin{equation}\label{eq:three-dimensional-local}
    \|(H_{\mathrm{cyl}}+V_b^{\mathrm{cap}}-\lambda)Q_b(\lambda)-\iota_b\|
\le Z^{-150},\qquad
  \supp Q_b(\lambda)f\subset B_{9Z/4}.
      \end{equation}
    The same excluded energies $\mathcal Z_b$ work for every $f\in L^2(B_b)$. The potential does not
  depend on $E$.
   \end{proposition}

  \begin{proof}
    We set $D_S\coloneqq B_{3S}\subset\R^2$, $I_b\coloneqq(-Z-w_b,Z+w_b)$, and $\Omega_b\coloneqq D_S\times I_b$. We define the
	Dirichlet operators
   \[
   A_S\coloneqq-\Delta_\rho+U,\qquad A_S^0\coloneqq-\Delta_\rho+V_{b_n}
\quad\text{on }D_S,\qquad L_b\coloneqq-\partial_z^2+v_b\quad\text{on }I_b
  \]
      by their nonnegative quadratic forms. On $\Omega_b$, we set
    \[
\mathscr H_b\coloneqq A_S\otimes I+I\otimes L_b,\qquad
	 \mathscr H_b^0\coloneqq A_S^0\otimes I+I\otimes L_b.
      \]
      Since $A_b=1$ on $D_S$, the operator $\mathscr H_b$ is the Dirichlet restriction of $H_{\mathrm{cyl}}+V_b^{\mathrm{cap}}$.

  We apply Proposition~\ref{prop:planar-block} on
 $B_{b_n}$ with $B\coloneqq640$ and $T=6$. Its excluded
   set $Y_n\subset[S^{-1/4},E]$ satisfies $|Y_n|\le
      S^{-6}$ at all sufficiently large stages. Indeed,
$b_n\asymp S^{1/3}\ge cb^{10/3}>b$, so $K(b_n)\ge
      K(b)\ge E$ eventually. The threshold for the
    actual ceiling $K(b_n)$ supplies these bounds at the transverse
 source radius $b_n$.

   For the axial inverse, we use Lemma~\ref{lem:axial-inverse} with $L=8$ and its actual ceiling $k_b\coloneqq K(b)$, and restrict its exceptional set to $[0,E]$. We
denote this set by $\mathcal E_b^{\mathrm{ax}}$ and its constants by $C_{k_b},c_{k_b}$. The threshold choice in \eqref{eq:diagonal-axial-bounds} ensures
	\begin{equation}\label{eq:cylinder-ceiling-bounds}
C_{k_b}\le b,\qquad
  C_{k_b}Z^{12}e^{-c_{k_b}w_b}\le Z^{-152}.
\end{equation}
The constants can include the finitely many localized derivative bounds at that
 ceiling. Constants depending only on the fixed bound $E$ will instead be denoted
    by $C_E$.

	We write $(\varepsilon_j,\psi_j)_{j\ge1}$ for a complete orthonormal
      eigenbasis of $L_b$, and $(E_k,\phi_k)_{k\ge1}$ for one of
$A_S$, including multiplicities. We exclude from $[0,E]$ the $Z^{-4}$-neighborhoods of both product spectra, all translates $\varepsilon_j+Y_n$ with $\varepsilon_j\le E$, and all translates $E_k+\mathcal E_b^{\mathrm{ax}}$
    with $E_k\le E$. This defines $\mathcal Z_b$. We use this same set for every $f\in L^2(B_b)$.

      For a compact-resolvent operator $T$, we write $N_T(E)\coloneqq\rank\1_{(-\infty,E]}(T)$. The min--max principle gives
 \[
     \begin{gathered}
  N_{A_S}(E)+N_{A_S^0}(E)\le C_ES^2,\qquad N_{L_b}(E)\le C_EZ,\\
 N_{\mathscr H_b}(E)+N_{\mathscr H_b^0}(E)
\le\bigl(N_{A_S}(E)+N_{A_S^0}(E)\bigr)N_{L_b}(E)\le C_ES^2Z.
  \end{gathered}
	\]
 The same bounds hold at $E+1$, as needed for the spectral neighborhoods. Consequently
\[
\begin{aligned}
 |\mathcal Z_b| \le C_E\bigl(S^2Z^{-3}+ZS^{-6}+C_{k_b}S^2Z^{-7}\bigr) \le C_E\bigl(b^{-60}+b^{-20}+b^{-219}\bigr)\le C_Eb^{-1}.
 \end{aligned}
\]
Here we used \eqref{eq:cylinder-ceiling-bounds} in the third term. For $\lambda\in[0,E]\setminus\mathcal Z_b$, both product resolvents have norm at most $Z^4$. For either planar operator $A\in\{A_S,A_S^0\}$,
      \[
	 \sigma(A)+\varepsilon_j\subset\sigma(A\otimes I+I\otimes L_b)
      \ \Longrightarrow\ \dist(\lambda-\varepsilon_j,\sigma(A))\ge Z^{-4}.
    \]
 We set $\mu_j\coloneqq\lambda-\varepsilon_j$. The same
	bound holds for every $j$, even when $\mu_j$ is
arbitrarily negative, so
	\begin{equation}\label{eq:shifted-inverses}
	\|(A_S-\mu_j)^{-1}\|+\|(A_S^0-\mu_j)^{-1}\|\le2Z^4
     \quad\text{for }j\ge1.
\end{equation}
	We extend $f\in L^2(B_b)$ by zero to $\Omega_b$ and define $f_j(\rho)\coloneqq\int_{I_b}\overline{\psi_j(z)}f(\rho,z)\,\dd z$ and $u_j\coloneqq(A_S-\mu_j)^{-1}f_j$. We have $f_j(\rho)=0$ for $|\rho|\ge b$. The orthonormality relation is $\int_{I_b}\overline{\psi_j(z)}\psi_k(z)\,\dd z=\delta_{jk}$. Since
 $(\psi_j)_{j\ge1}$ is complete, Parseval and integration in $\rho$ give
\[
 \|f\|_{L^2(\Omega_b)}^2
     =\int_{D_S}\sum_{j\ge1}|f_j(\rho)|^2\,\dd\rho
  =\sum_{j\ge1}\|f_j\|_{L^2(D_S)}^2.
   \]
 The product decomposition now gives
      \[
   u\coloneqq(\mathscr H_b-\lambda)^{-1}f
=\sum_{j\ge1}u_j(\rho)\psi_j(z),\qquad
\|u\|\le Z^4\|f\|.
      \]
      The series converges in $L^2(\Omega_b)$. In particular, it includes all axial eigenfunctions, including those for which
      $\mu_j$ is negative.

   We first estimate $u$ near the boundary $|\rho|=3S$. We set $\mathcal O_S\coloneqq\{9S/4<|\rho|<3S\}$. If $\mu_j\ge S^{-1/4}$, it lies outside $Y_n$. For $g\in L^2(B_{b_n})$, extended by zero to $D_S$, the planar inverse
	gives
\[
\begin{aligned}
q^0&\coloneqq Q_{b_n}(\mu_j)g,\quad r^0\coloneqq(A_S^0-\mu_j)q^0-g,
  \quad\|r^0\|\le S^{-B}\|g\|,\\
  (A_S^0-\mu_j)^{-1}g&=q^0-(A_S^0-\mu_j)^{-1}r^0,
  \qquad\1_{\mathcal O_S}q^0=0.
\end{aligned}
   \]
      Thus \eqref{eq:shifted-inverses} gives $\|\1_{\mathcal O_S}(A_S^0-\mu_j)^{-1}g\|\le Z^4S^{-B}\|g\|$. The preceding potentials
  form a bounded sum $U_{<n}\coloneqq\sum_{m<n}V_{b_m}$ supported in $B_{b_n}$. The next annular potential lies outside $D_S$, so
 $A_S=A_S^0+U_{<n}$. Also $b_n\asymp S^{1/3}\ge
    cb^{10/3}>b$. The identity
	\[
(A_S-\mu_j)^{-1}f_j
     =(A_S^0-\mu_j)^{-1}
  \bigl(f_j-U_{<n}(A_S-\mu_j)^{-1}f_j\bigr)
\]
 therefore gives
      \begin{equation*}
\|\1_{\mathcal O_S}u_j\|
 \le Z^4S^{-B}\bigl(\|f_j\|+\|U_{<n}\|\|u_j\|\bigr)
   \le Z^4S^{-B}(1+CZ^4)\|f_j\|
\le CZ^8S^{-B}\|f_j\|.
\end{equation*}

For the remaining values of $\mu_j$, we use a smaller disk that contains $B_b$ and has low-energy
 threshold above $S^{-1/4}$. We set $S_0\coloneqq R_{n-2}\asymp S^{1/9}$ and note that
   \[
 B_b\subset B_{S_0/2},\qquad \mu_j<S^{-1/4}\le cS_0^{-2}
      \quad\text{for sufficiently large }b,
\]
since
\[
S_0\ge cb^{10/9},\qquad \frac b{S_0}\le Cb^{-1/9}\to0,\qquad
  S^{-1/4}S_0^2\le CS^{-1/36}\to0.
	\]
	Lemma~\ref{lem:low-inverse} therefore gives
      \[
      q_j\coloneqq Q_{S_0}^{\mathrm{in}}(\mu_j)f_j,\qquad
	r_j\coloneqq(h_0-\mu_j)q_j-f_j,\qquad
\|r_j\|\le Ce^{-cS_0^{1/6}}\|f_j\|.
\]
  Thus $u_j=q_j-(A_S-\mu_j)^{-1}r_j$ and $\1_{\mathcal O_S}q_j=0$. By \eqref{eq:shifted-inverses},
\begin{equation*}
\|\1_{\mathcal O_S}u_j\|
     \le CZ^4 e^{-cS_0^{1/6}}\|f_j\|.
\end{equation*}
      We combine the two estimates by setting
    \[
   d_b\coloneqq Z^8S^{-B}+Z^4e^{-cS_0^{1/6}}.
\]
They give $\|u_j\|_{L^2(\mathcal O_S)}\le Cd_b\|f_j\|$ for every $j$. To sum these bounds, we first take a positive integer $K$ and
set $u^{(K)}\coloneqq\sum_{j=1}^Ku_j\psi_j$. Orthonormality gives
   \[
\int_{\mathcal O_S}\int_{I_b}|u^{(K)}(\rho,z)|^2\,\dd z\,\dd\rho
=\sum_{j,k=1}^K\int_{\mathcal O_S}
      \overline{u_j(\rho)}u_k(\rho)\,\dd\rho
\int_{I_b}\overline{\psi_j(z)}\psi_k(z)\,\dd z
=\sum_{j=1}^K\|u_j\|_{L^2(\mathcal O_S)}^2.
     \]
We let $K\to\infty$. The left side converges because $u^{(K)}\to u$ in $L^2(\Omega_b)$. Thus
  \[
    \|u\|_{L^2(\mathcal O_S\times I_b)}^2
    =\sum_{j\ge1}\|u_j\|_{L^2(\mathcal O_S)}^2
  \le Cd_b^2\sum_{j\ge1}\|f_j\|^2=Cd_b^2\|f\|^2.
      \]
    We next estimate the derivatives needed to cut off $u$ in $\rho$. We choose a
	radial cutoff $\Theta$, equal to one on $B_{5S/2}$ and zero outside $B_{11S/4}$,
	with derivatives of order $k\le2$ bounded by $C_kS^{-k}$. We also choose $\chi\in
    C_c^\infty(\mathcal O_S)$, equal to one near
      $\supp\nabla\Theta\cup\supp\Delta\Theta$, with $|\nabla\chi|\le C/S$. On $\mathcal O_S$ the equation for $u_j$ is
	$(-\Delta_\rho+U-\mu_j)u_j=0$. Testing it against $\chi^2u_j$ gives
   \[
\|\nabla_\rho(\chi u_j)\|^2
+\int_{D_S}U|\chi u_j|^2\,\dd\rho
	 -\mu_j\|\chi u_j\|^2
  =\|u_j\nabla\chi\|^2.
  \]
	Since $U\ge0$ and $\mu_j\le E$, this identity gives
  \[
\|\nabla_\rho u_j\|_{\supp\nabla\Theta}
    \le C\bigl(S^{-1}+\sqrt{\max\{\mu_j,0\}}\bigr)
\|u_j\|_{L^2(\mathcal O_S)}
\le C_{E}\|u_j\|_{L^2(\mathcal O_S)}.
\]
 This is also the localized energy estimate in
    \cite[Lemmas~3.5--3.6]{FLM}. The negative values of
  $\mu_j$ improve the left side of the tested identity,
      so the same constant works as $\mu_j\to-\infty$.
  We therefore have
      \[
     [-\Delta_\rho,\Theta]u_j
=-2\nabla\Theta\cdot\nabla_\rho u_j-(\Delta\Theta)u_j,
	 \qquad
      \|[-\Delta_\rho,\Theta]u_j\|
  \le CS^{-1}\|u_j\|_{L^2(\mathcal O_S)}.
	\]
  The cutoff depends only on $\rho$, so it does not change the axial functions $\psi_j$. Applying the same orthogonality calculation gives
  \[
\|[-\Delta_\rho,\Theta]u^{(K)}\|^2
  =\sum_{j=1}^K\|[-\Delta_\rho,\Theta]u_j\|^2
\le CS^{-2}d_b^2\sum_{j=1}^K\|f_j\|^2.
    \]
    The same estimate applied to a difference of partial sums shows that these commutators converge in $L^2$. Since $u^{(K)}\to u$ in $L^2$, their limit is $[-\Delta_\rho,\Theta]u$ in the distributional sense. We obtain
\begin{equation}\label{eq:transverse-commutator}
   \|[-\Delta_\rho,\Theta]u\|
    \le CS^{-1}d_b\|f\|
	 =CS^{-1}\bigl(Z^8S^{-B}+Z^4e^{-cS_0^{1/6}}\bigr)\|f\|.
  \end{equation}
The squared-norm identity gives a bound uniform in $K$, the number of axial terms.
   The $u_j(\rho)$ estimates hold on the full angular Fourier space, with no angular counting factor.

      We next estimate $u$ near the two ends of the cylinder. We use the other product decomposition,
 \[
 \begin{gathered}
f=\sum_{k\ge1}\phi_k(\rho)\widehat f_k(z),\qquad
      u=\sum_{k\ge1}\phi_k(\rho)\widehat u_k(z),\\
\nu_k\coloneqq\lambda-E_k,\qquad
     \widehat u_k=(L_b-\nu_k)^{-1}\widehat f_k,\qquad
    \sum_{k\ge1}\|\widehat f_k\|^2=\|f\|^2.
      \end{gathered}
      \]
      Every $\widehat f_k$ is supported in $(-b,b)\subset(-Z,Z)$. The product spectral exclusion gives $\|(L_b-\nu_k)^{-1}\|\le Z^4$ for every $k$. If $\nu_k\ge0$, then $E_k\le E$, and membership in $\mathcal
    E_b^{\mathrm{ax}}$ was excluded by the translate $E_k+\mathcal E_b^{\mathrm{ax}}$. If $\nu_k<0$, Lemma~\ref{lem:axial-inverse}
    requires no exclusion. Thus in every case we may set
    \[
    q_k\coloneqq Q_Z^{\mathrm{ax}}(\nu_k)\widehat f_k,\qquad
   r_k\coloneqq(L_b-\nu_k)q_k-\widehat f_k,
	 \qquad
   \|r_k\|\le C_{k_b}Z^8e^{-c_{k_b}w_b}\|\widehat f_k\|.
      \]
	We have $\widehat u_k=q_k-(L_b-\nu_k)^{-1}r_k$. Every $q_k$ vanishes on the same outer collar
   \[
   \mathcal T_b\coloneqq\{z:Z+w_b/2<|z|<Z+3w_b/4\}.
  \]
The resolvent bound and \eqref{eq:cylinder-ceiling-bounds} therefore give
\[
\|\widehat u_k\|_{L^2(\mathcal T_b)}
\le C_{k_b}Z^{12}e^{-c_{k_b}w_b}\|\widehat f_k\|
\le Z^{-152}\|\widehat f_k\|.
\]

 We choose a smooth even cutoff $\zeta$, equal to one
    for $|z|\le Z+9w_b/16$ and zero for $|z|\ge
Z+11w_b/16$. Its first two derivatives are supported strictly inside
	$\mathcal T_b$ and are bounded uniformly for large $b$. On this collar the equation is $(-\partial_z^2+v_b-\nu_k)\widehat u_k=0$.
    The one-dimensional localized energy estimate, using $v_b\ge0$ and $\nu_k\le E$, gives
   \[
\|\partial_z\widehat u_k\|_{\supp\zeta'}
	 \le C_E\|\widehat u_k\|_{L^2(\mathcal T_b)}.
    \]
	The negative values of $\nu_k$ improve the tested identity, so the constant is uniform for all $k$. The cutoff depends only on $z$. Applying the same orthogonality calculation as above, now to the complete planar eigenbasis, gives
    \begin{equation}\label{eq:axial-commutator}
	\|[-\partial_z^2,\zeta]u\|^2
=\sum_{k\ge1}\|[-\partial_z^2,\zeta]\widehat u_k\|^2
  \le C_EZ^{-304}\sum_{k\ge1}\|\widehat f_k\|^2
	 =C_EZ^{-304}\|f\|^2.
    \end{equation}
The bound gives $L^2$ convergence of the finite-sum commutators, hence the identity for all planar eigenvalues,
   including the infinitely many giving negative $\nu_k$. No counting factor occurs.

We define $Q_b(\lambda)f\coloneqq\Theta(\rho)\zeta(z)u$ and extend it by zero. Both cutoffs equal one on
      $B_b$. The residual has just two terms:
    \[
(H_{\mathrm{cyl}}+V_b^{\mathrm{cap}}-\lambda)Q_bf-f
=\zeta[-\Delta_\rho,\Theta]u
+\Theta[-\partial_z^2,\zeta]u.
   \]
      The function $Q_bf$ is compactly supported in $\Omega_b$, where the two differential expressions agree. Its distributional
   Laplacian belongs to $L^2$. The Fourier transform then gives $Q_bf\in H^2(\R^3)$, with support in $B_{9Z/4}$. Now
\eqref{eq:transverse-commutator}--\eqref{eq:axial-commutator} give
   \[
\|(H_{\mathrm{cyl}}+V_b^{\mathrm{cap}}-\lambda)Q_b-\iota_b\|
	\le C_E\bigl(Z^8S^{-641}
+Z^4S^{-1}e^{-cS_0^{1/6}}+Z^{-152}\bigr).
      \]
  The polynomial error is smaller than $Z^{-150}$,
 since
\[
     Z^{158}S^{-641}
	\le b^{40\cdot158-10\cdot641}=b^{-90}\longrightarrow0.
	\]
      Since $S_0\ge cb^{10/9}$, the remaining terms satisfy
    \[
Z^{154}S^{-1}e^{-cS_0^{1/6}}\longrightarrow0,
\qquad Z^{150}Z^{-152}\longrightarrow0.
\]
  This proves \eqref{eq:three-dimensional-local} for sufficiently large $b$.
    \end{proof}

\subsection{The three-dimensional potential}

	We now sum the pairs of disk barriers and apply the criterion from Section~\ref{sec:spectral-averaging}. Since
    \eqref{eq:cylinder-scales} gives $S\le CZ^{3/4}$, the cost of one pair satisfies
\begin{equation}\label{eq:cap-cost}
     \mathcal I_{3,\sigma}(V_b^{\mathrm{cap}})
 \le9Z^{\sigma-2}|\supp V_b^{\mathrm{cap}}|
     \le C(\log Z)^2S^2Z^{\sigma-2}
     \le C(\log Z)^2Z^{\sigma-1/2}.
    \end{equation}
 For $0\le\sigma<1$, we set $c_\sigma\coloneqq\int_\R(1+t^2)^{(\sigma-2)/2}\,\dd t<\infty$. We integrate in $z$ to obtain
\begin{equation*}
\mathcal I_{3,\sigma}(U(\rho))=c_\sigma\mathcal I_{2,\sigma}(U).
     \end{equation*}

Starting from a large $\widetilde b_1$, we set $Z_j\coloneqq\widetilde b_j^{40}$ and $\widetilde
  b_{j+1}\coloneqq3Z_j$, and add the corresponding pairs of disk barriers. Their
supports are disjoint in $z$. By \eqref{eq:cap-cost}, their weighted integrals are
 summable for $0\le\sigma<1/3$. We choose the exponents in
Proposition~\ref{prop:assembly} once, independently of the upper energy endpoint:
   \[
d=3,\qquad p\coloneqq40,\qquad \alpha\coloneqq\tfrac72,\qquad
  \beta\coloneqq150,\qquad Q_3\coloneqq6000.
   \]
These choices satisfy $\beta+1=151>(p+2)\alpha=147$ and $Q_3=6000>p^2\alpha=5600$. For every fixed $E>0$, Proposition~\ref{prop:cylinder-block} applies. The measures of the exceptional sets are summable. For $t>0$, we define
   \[
   W_3(x)\coloneqq\langle x\rangle^{-2Q_3},\qquad
   V_3(\rho,z)\coloneqq U(\rho)
   +\sum_jV_{\widetilde b_j}^{\mathrm{cap}}(\rho,z)+tW_3(\rho,z).
    \]
      For every fixed $E>0$, the estimates hold at all
sufficiently large stages. We take the exceptional set
   to be the entire interval at each earlier stage.
   Proposition~\ref{prop:assembly} then gives pure
   point spectral type on $[0,E]$ for almost every
	$t>0$. We intersect these full-measure sets along the sequence $E=1,2,\ldots$. For
every $t$ in the resulting full-measure set, $H_{V_3}$ has pure point spectral type
   on $(0,\infty)$, with finite eigenvalue multiplicities. We choose $t>0$ in this
set sufficiently small when estimating the weighted integrals below. Since the disk
      barriers have disjoint supports, Minkowski's inequality gives, also at
   $\sigma=0$,
    \[
 \sqrt{\mathcal I_{3,\sigma}(V_3)}
   \le\sqrt{c_\sigma\mathcal I_{2,\sigma}(U)}
+\sqrt{\sum_j\mathcal I_{3,\sigma}(V_{\widetilde b_j}^{\mathrm{cap}})}
   +t\sqrt{\mathcal I_{3,\sigma}(W_3)}.
 \]
We make both integrals arbitrarily small by choosing large initial radii and then small $t$.

   Nonnegativity gives $\sigma(H_{V_3})\subset[0,\infty)$. For the reverse inclusion, we set
	\[
    \rho_k\coloneqq5R_k\mathbf e_1,\qquad
     \ell_k\coloneqq R_k^{1/2},\qquad
	a_k\coloneqq(\rho_k,0),\qquad \mathcal B_k\coloneqq B_{\ell_k}(a_k).
	\]
      For $(\rho,z)\in\mathcal B_k$, we have
    \[
 5R_k-\ell_k<|\rho|<5R_k+\ell_k,\qquad |z|<\ell_k.
   \]
   Since $R_{k+1}\ge10R_k$, this planar annulus lies between the $k$th and $(k+1)$st annular potentials
 for large $k$. Hence $U(\rho)=0$ throughout $\mathcal B_k$. We write $S_j$ for the planar radius used in
    the $j$th pair of disk barriers. Their defining scales give
 \[
	 S_j\le C\widetilde b_j^{30}=CZ_j^{3/4},\qquad
     \supp V_{\widetilde b_j}^{\mathrm{cap}}
   \subset\{|\rho|\le4S_j,\ |z|>Z_j\}.
  \]
We compare this support with $\mathcal B_k$ in two cases. If $Z_j>\ell_k$, then every point of the barrier has $|z|>\ell_k$, whereas every point of $\mathcal B_k$ has
$|z|<\ell_k$. If $Z_j\le\ell_k$, then every point of the barrier has
    \[
|\rho|\le CZ_j^{3/4}\le C\ell_k^{3/4}
=CR_k^{3/8}<5R_k-\ell_k
  \]
for all sufficiently large $k$, uniformly in $j$. Every point of
  $\mathcal B_k$ has $|\rho|>5R_k-\ell_k$. Thus no disk barrier meets
  $\mathcal B_k$. The remaining potential satisfies
\[
\sup_{x\in\mathcal B_k}t\langle x\rangle^{-2Q_3}
     \le t\langle5R_k-\ell_k\rangle^{-2Q_3}\longrightarrow0.
	 \]

      For $\lambda\ge0$, we choose $\xi\in\R^3$ with $|\xi|^2=\lambda$ and $\chi\in C_c^\infty(B_1)$ with
 $\|\chi\|_2=1$. Then
   \[
\begin{gathered}
	u_k(x)\coloneqq\ell_k^{-3/2}\chi\bigl((x-a_k)/\ell_k\bigr)e^{i\xi\cdot x},
    \qquad\|u_k\|_2=1,\qquad u_k\rightharpoonup0,\\
\|(H_{V_3}-\lambda)u_k\|_2
     \le\frac{2|\xi|}{\ell_k}\|\nabla\chi\|_2
+\frac1{\ell_k^2}\|\Delta\chi\|_2
  +\|V_3\|_{L^\infty(B_{\ell_k}(a_k))}\longrightarrow0.
    \end{gathered}
\]
      Hence $\sigma(H_{V_3})=\sigma_{\mathrm{ess}}(H_{V_3})=[0,\infty)$.

\subsection{Tensor sums}

      For $d=2r$, we use $r\ge1$ planar factors. For $d=3+2r$, we use one three-dimensional factor and $r\ge0$ planar factors. We write $k$ for the number of factors, with dimensions $m_1,\ldots,m_k\in\{2,3\}$, $\sum_jm_j=d$, and coordinates $x_j\in\R^{m_j}$. Set
      \[
V_d(x_1,\ldots,x_k)\coloneqq\sum_{j=1}^k V_{m_j}(x_j),\qquad
  H_d\coloneqq\sum_{j=1}^k I\otimes H_{V_{m_j}}\otimes I.
    \]
The operator identity holds on $H^2$, since the potentials are bounded. For $m\in\{2,3\}$, $d>m$, and
    $0\le\sigma<m-1$, we define $c_{d,m,\sigma}$ below. We substitute $z=|y|w$ for $y\in\R^m\setminus\{0\}$ and apply Minkowski to obtain
 \[
  \begin{gathered}
  c_{d,m,\sigma}\coloneqq\int_{\R^{d-m}}(1+|z|^2)^{(1-d+\sigma)/2}\,\dd z<\infty,
     \qquad c_{d,d,\sigma}\coloneqq1,\\
   \int_{\R^{d-m}}(|y|^2+|z|^2)^{(1-d+\sigma)/2}\,\dd z
   =c_{d,m,\sigma}|y|^{1-m+\sigma}, \qquad
	 \mathcal I_{d,\sigma}(V_d)^{1/2}
\le\sum_j\bigl(c_{d,m_j,\sigma}\mathcal I_{m_j,\sigma}(V_{m_j})\bigr)^{1/2}.
  \end{gathered}
   \]
    We apply this bound at $\sigma$ and at zero. Choosing the factors with
      sufficiently small weighted integrals makes both integrals arbitrarily small.

For either dimension $m\in\{2,3\}$, positivity gives $H_{V_m}u=0$ which shows $u=0$. Thus each factor is nonnegative, has zero kernel, and has pure point spectral type on $(0,\infty)$. Its normalized eigenfunctions form a
  complete orthonormal basis. The tensor products of these basis functions form a complete orthonormal basis of eigenfunctions of $H_d$, and their eigenvalues are positive sums of factor eigenvalues. Consequently $P_c(H_d)=0$
and $\ker H_d=\{0\}$. We also obtain
\[
\sigma(H_d)=\overline{\sigma(H_{V_{m_1}})+\cdots+\sigma(H_{V_{m_k}})}
   =[0,\infty)=\sigma_{\mathrm{ess}}(H_d),
\]
   since this half-line has no isolated spectral points. It remains
to make the supremum norm small at the same time as the two weighted
integrals. For $a>0$, set
   \[
V^{(a)}(x)\coloneqq aV_d(\sqrt a\,x),\qquad
(\mathcal U_af)(x)\coloneqq a^{d/4}f(\sqrt a\,x).
	 \]
   The map $\mathcal U_a$ is unitary, and direct calculation gives
\[
     \begin{gathered}
     H_{V^{(a)}}=a\mathcal U_aH_d\mathcal U_a^{-1},\qquad
	\|V^{(a)}\|_\infty=a\|V_d\|_\infty, \qquad 
\mathcal I_{d,\tau}(V^{(a)})
 =a^{(3-\tau)/2}\mathcal I_{d,\tau}(V_d),\qquad \tau\in\{0,\sigma\}.
      \end{gathered}
	 \]
All three terms in \eqref{eq:main-weight} tend to zero as $a\downarrow0$. The unitary identity preserves pure point spectral type on
    the positive half-line, the zero kernel, and the spectrum $[0,\infty)$. In dimensions two and three it also preserves the finite
eigenvalue multiplicities furnished by Proposition~\ref{prop:assembly}. Taking $a$ sufficiently small proves the theorem.

\par
  \end{document}